\documentclass[oneside]{amsart}
\usepackage[T1]{fontenc}
\usepackage[utf8]{inputenc}
\usepackage{color}
\usepackage{mathtools}
\usepackage{amstext}
\usepackage{amsthm}
\usepackage{verbatim}
\usepackage{amssymb}
\usepackage{esint}
\usepackage[unicode=true,
 bookmarks=false,
 breaklinks=false,pdfborder={0 0 1},colorlinks=false]
 {hyperref} 
\usepackage{a4wide}
\usepackage{appendix}

\makeatletter
\numberwithin{equation}{section}
\numberwithin{figure}{section}

\@ifundefined{date}{}{\date{}}

\usepackage{mathrsfs}
\usepackage{graphicx}
\usepackage{enumitem}
\usepackage{tikz}
\usepackage{color}

\theoremstyle{definition}
\newtheorem{definition}{Definition}[section]
\newtheorem{cor}[definition]{Corollary}

\theoremstyle{plain}
\newtheorem{thm}[definition]{Theorem}\theoremstyle{plain}
\newtheorem{lem}[definition]{Lemma}\theoremstyle{remark}
\newtheorem{rem}[definition]{Remark}\theoremstyle{plain}
\newtheorem{prop}[definition]{Proposition}
\newtheorem{mainthm}{Theorem}

\numberwithin{equation}{section}

\newcommand{\Zb}{\mathbb{Z}}
\newcommand{\Rb}{\mathbb{R}}

\newcommand{\Nb}{\mathbb{N}}

\newcommand{\Eb}{\mathbb{E}}

\newcommand{\Ab}{\mathbb{A}}
\newcommand{\Db}{\mathbb{D}}

\newcommand{\FS}{\mathcal{S}}

\newcommand{\D}{\mathcal{D}}

\newcommand{\N}{\mathcal{N}}

\renewcommand{\S}{\mathcal{S}}
\newcommand{\F}{\mathcal{F}}

\newcommand{\Pc}{\mathcal{P}}
\newcommand{\T}{\mathcal{T}}
\newcommand{\R}{\mathcal{R}}

\newcommand{\norm}[1]{\left\| #1 \right\|}
\newcommand{\md}[1]{\left| #1 \right|}

\newcommand{\p}[1]{\left( #1 \right)}
\newcommand{\set}[1]{\left\lbrace #1 \right\rbrace}

\newcommand{\les}{\lesssim}

\newcommand{\e}{\varepsilon}

\newcommand{\Ind}{\mathbf{1}}

\newcommand{\Int}{\mathrm{Int}}
\newcommand{\Cl}{\mathrm{Cl}}
\newcommand{\Adm}{\mathrm{Adm}}

\newcommand{\Tl}{\overline{T}}

\newcommand{\Tcl}{\overline{\mathcal{T}}}

\newcommand{\bd}[1]{\underline{\mathrm B} (#1)}

\title{Noncommutative maximal ergodic inequalities for amenable groups }

\author{Léonard Cadilhac}

\address{Institut de Mathématiques de Jussieu, Sorbonne Université, 4 place de Jussieu, 75252 Paris Cedex 05, France}

\email{cadilhac@imj-prg.fr}

\author{Simeng Wang}

\address{Institute for Advanced Study in Mathematics, Harbin Institute of Technology, Harbin 150001, China}

\email{simeng.wang@hit.edu.cn}

\makeatother

\begin{document}

\begin{abstract}
We prove a pointwise ergodic theorem and a maximal inequality for actions of amenable groups on noncommutative measure spaces. To do so, we establish a square function estimate quantifying the difference between ergodic averages and some conditional expectations. 

Our main technical results are the construction of a well-behaved filtration, based on the quasi-tilings of Ornstein and Weiss, and the square function bound, which we derive from non-doubling noncommutative Calder\'on-Zygmund decomposition. 

For actions on usual measure spaces, we obtain new variational ergodic inequalities and jump estimates.
\end{abstract}

\maketitle
\global\long\def\md#1{\left|#1\right|}
 \global\long\def\Pc{\mathcal{P}}
 \global\long\def\M{\mathcal{M}}
 \global\long\def\norm#1{\|#1\|}
 \global\long\def\les{\lesssim}
 \global\long\def\e{\mathcal{\varepsilon}}
 \global\long\def\p#1{\left(#1\right)}
 \global\long\def\N{\mathcal{N}}
 \global\long\def\Ind{\mathcal{\mathbf{1}}}
 \global\long\def\Eb{\mathcal{\mathbf{\mathbb{E}}}}
 \global\long\def\F{\mathcal{\mathbf{\mathcal{F}}}}

\section{Introduction}

\label{sec:intro}

Pointwise ergodic theorems play a fundamental role in ergodic theory. In particular, the celebrated Birkhoff ergodic theorem establishes the pointwise convergence of ergodic averages for a measure-preserving transformation. More precisely, let $T$ be a measure-preserving transformation of a measure space $(X,\mu)$. For a function $f$ on $X$ and $n\in \Nb$, we consider the ergodic averages
$$
A_n(f) := \dfrac1{n}\sum_{k=0}^{n-1} f\circ T^k.
$$
In 1931, Birkhoff proved the pointwise ergodic theorem showing that $\Ab_n(f)$ converges almost everywhere for all $f\in L_p(X,\mu)$ with $1\leq p < \infty$. When multiple transformations are involved, the study of ergodic theorems for group actions naturally arises. Let $G$ be a locally compact group, $r$ a right-invariant Haar measure on $G$, and $\alpha$ a continuous action of $G$ on $X$ by $\mu$-preserving transformations. In this context, there is often no canonical way to define ergodic averages; however, given a sequence $(F_n)_{n\in \Nb}$ of compact subsets of $G$, we may set 
$$
A_{n}(f) := \dfrac1{r(F_n)}\int_{F_n} f \circ \alpha_g dr(g) 
$$
and study the corresponding pointwise convergence.
In this line of investigation, the first major step was taken by Calderón \cite{Cal53}, who considered averages based on sequences $(F_n)_{n\in\Nb}$ that are  {\it asymptotically invariant}, a condition that requires $G$ to be amenable.  Since then, amenability has remained a key property in understanding pointwise ergodic theorems for group actions, although techniques for certain non-amenable groups have since emerged --- see \cite{Nev06} for a survey of the state of the art in 2006. In particular, in 2001, Lindenstrauss \cite{Lin01} solved a longstanding problem, showing that every amenable group admits  a pointwise ergodic sequence, {\it i.e.}, a sequence $(F_n)_{n\in \Nb}$ such that $(A_n(f))$ converges almost everywhere to a $G$-invariant function for any $f \in L_p(X,\mu)$ with $ 1\leq p < \infty$. 

The purpose of this paper is to explore analogous problems within the context of actions on noncommutative measure spaces, that is, on von Neumann algebras. Ergodic theory and von Neumann algebras have had a profound interaction since the beginning of the theory of operator algebras. The study of pointwise and maximal ergodic theorems in the noncommutative setting was initiated in the 1970s in Lance's seminal work \cite{Lan76}, and his result was subsequently improved by Conze, Dang-Ngoc, Kummerer, Yeadon and others \cite{ConzeDangNgoc1978,Kummerer1978,Yeadon1977}. In 2007, Junge and Xu \cite{JX07} established the Birkhoff (or more generally, Dunford-Schwartz) maximal ergodic theorems on noncommutative $L_p$-spaces, thereby resolving a longstanding open problem in the field of noncommutative ergodic theory. This breakthrough led to further progress (see e.g. \cite{Bekjan2008,HongSun2018,CGP20max,HLX24single}); in particular, the noncommutative Stein ergodic theorem for free group actions was investigated in \cite{AnantharamanDelaroche2006,Hu2008}. Since then, the research for relevant results for actions of amenable groups has become a natural open problem in the field. Classical approaches do not extend directly to the noncommutative setting due to the current lack of understanding of noncommutative covering lemmas.  The first general results in this direction were presented by the second author, Hong and Liao in \cite{HLW21}, where a noncommutative analog of the aforementioned Calderón ergodic theorem \cite{Cal53} was established. The work \cite{HLW21} relies essentially on the doubling condition of the underlying groups, and its methods do not apply to the study of the noncommutative Lindentrauss ergodic theorem for general actions of amenable groups.

In this paper, we establish a noncommutative pointwise ergodic theorem for actions by general amenable groups. We begin by setting up the framework in the noncommutative setting. Consider a von Neumann algebra $\M$ equipped with a trace $\tau$, which is understood as a noncommutative measure space. In this context, measure-preserving transformations are naturally replaced by trace preserving $\ast$-automorphisms of $\M$, and the almost everywhere convergence is replaced by the so-called {\it (bilaterally) almost uniform convergence} (see Definition~\ref{def:au}), which were introduced in Lance's aforementioned pioneering work \cite{Lan76}. We will be interested in weakly continuous actions $\alpha\coloneqq (\alpha_g )_{g\in G}$  of $G $ on $ \M$ by $\tau $-preserving $\ast$-automorphisms. For $x\in \mathcal{M}$ and a given suitable sequence  $(F_n)_{n\geq 0}$ of compact subsets in $G$, the ergodic averages are now defined by
$$
A_{n}(x) := \dfrac1{r(F_n)}\int_{F_n} \alpha_g(x) dr(g).
$$
We extend Lindenstrauss's pointwise theorem (\cite{Lin01}) in the following way. 

\begin{thm}\label{thm:intromain}
	Let $G$ be amenable and $\alpha$ be as above. There is a F\o lner sequence $(F_n)_{n\geq 0}$ in $G$ such that for any $x\in L_{p}(\mathcal{M})$ with  $1\leq p<2$ (resp. $2\leq p < \infty$),  	\[
	A_{n}(x) \rightarrow Px  \quad \text{bilateral almost uniformly (resp. almost uniformly)},
	\]
where $P$ is the projection onto the space of fixed points of $\alpha$.
\end{thm}
We refer to Section \ref{sec:prelim} for details and a more complete statement (Theorem \ref{thm:main}). Note that Lindenstrauss's original theorem \cite{Lin01} applies to \emph{tempered} F\o lner sequences. Our result is similar in  that it guarantees the existence of a pointwise ergodic sequence $(F_n)_{n \in \mathbb{N}}$, but requires the sequence to be \emph{admissible}---a more restrictive and technical condition (see Definition~\ref{def:filtered}). However, this approach allows us to establish variation and jump estimates  for $(A_n)_{n \in \mathbb{N}}$ that provide quantitative refinements of pointwise convergence and strengthen Lindenstrauss's maximal inequality.
These estimates are new even in the classical setting 
(see Theorem~\ref{thm:quantitative}).

The proofs rely on several new geometric and analytic methods. As in \cite{HLW21} and \cite{Lin01}, the core problem is to establish an analogue of the Hardy--Littlewood maximal inequality for operator-valued functions on amenable groups. In the absence of covering lemmas, a well-tested method is to compare ergodic averages to martingales. To this end, we take inspiration from the method developed by Jones, Kaufman, 
Rosenblatt, and Wierdl~\cite{JKRM88osc,JRW03} 
for oscillation and variation estimates, which further dates back to Bourgain \cite{Bou89}.  Hong and Xu \cite{HX21,Xu21} have shown that this method 
can be extended to the operator-valued functions on Euclidean spaces. However, their method relies essentially on geometric arguments in $\Rb^d$: it considers ball averages over $F_n = B(0,2^n)$ and exploits their relationship with the dyadic filtration $\Pc_n := \set{[0,2^n)^d + 2^nk : k\in \Zb^d}$, where the metric structure and doubling geometry play a crucial role in their estimates.

The main achievement of this paper is to develop a new geometric framework 
for noncommutative Calder\'on--Zygmund theory 
and dyadic-like martingale methods on abstract amenable groups, 
in which the aforementioned key ideas from the Euclidean setting---such as 
the interplay between ergodic averages and martingales---are realized through intrinsic group-theoretic constructions, 
without relying on any ambient metric structure.
This framework draws on techniques from several research directions:

\begin{itemize}
    \item The ideas of Jones, Kaufman, Rosenblatt, Mirek, and Wierdl 
          \cite{JKRM88osc,JRW03} on variational inequalities, 
          which introduce a square function relating ergodic averages 
          and conditional expectations. 
          Their approach however relies on Euclidean geometric arguments 
          and Calder\'on--Zygmund decomposition, which requires additional effort to extend to the setting of abstract groups;

    \item Non-doubling Calder\'on--Zygmund theory, 
          popularized by Tolsa \cite{Tol01,Tol03}, 
          motivated by the study of non-doubling measures on $\mathbb{R}^d$. 
          We adapt this theory to the setting of Haar measures 
          on locally compact amenable groups;

    \item Quasi-tilings of amenable groups, 
          introduced by Ornstein and Weiss \cite{OW87}, 
          as a key tool in understanding the ergodic properties 
          of amenable group actions. 
          To the best of our knowledge, 
      this is the first time quasi-tilings are applied 
      to maximal inequalities; 
      we use them here to construct a proxy for the dyadic filtration;

    \item Noncommutative analogs of the first two theories, 
          particularly the non-doubling noncommutative Calder\'on--Zygmund 
          decomposition developed in \cite{CCP21}, 
          as well as the operator-valued extensions of variational ergodic inequalities 
          in \cite{HX21,Xu21}.
\end{itemize}

The remainder of this paper is organized as follows. In Section \ref{sec:prelim}, we first introduce the necessary preliminaries and outline the structure of the proof of Theorem \ref{thm:intromain}. The detailed discussions of key technical components are postponed to subsequent sections. In Section \ref{sec:geometry}, we construct admissible F{\o}lner sequences for amenable groups. Section \ref{sec:CZdec} presents the noncommutative Calderón–Zygmund decomposition in the abstract group-theoretic setting. Finally, in Section \ref{sec:diffop}, we establish the critical estimates for square functions.

\section{Preliminaries and scheme of the proof}
\label{sec:prelim}
We will first present in this section the entire scheme of the proof of Theorem \ref{thm:intromain}, and leave the detailed proof of several key results to later sections. To this end, we also include preliminaries on noncommutative integration and amenability.
\subsection{Noncommutative integration}
\label{sub:prelim_nc}

In this paper, $\M$ denotes a semifinite von Neumann algebra and $\tau$ a normal semifinite faithful trace on $\M$.
For $p\in (0,\infty)$, we consider the $L_p$-norm
$$\|x\|_p = [\tau (|x|^p )]^{1/p},\quad x\in \mathcal M, $$
and the associated noncommutative $L_p$-space $L_p(\M)$ on $\mathcal M$ is defined to be the completion of ${\{x\in \mathcal M : \|x\|_p <\infty\}}$ with respect to the norm $\|\,\|_p$. For convenience, we write $L_\infty (\mathcal M) =\mathcal M$. The reader is referred to \cite{PX03survey} for more information on noncommutative $L_p$-spaces. We will also consider the quasi norm given by
$$\|x\|_{p,\infty} = \sup_{t>0} t [\tau(\mathbf{1}_{[0,t]} (x))]^{1/p},\quad x\in \mathcal M, $$
and denote by $L_{p,\infty} (\mathcal M)$ the noncommutative weak $L_p$-space obtained by the completion of $\{x\in \mathcal M : \|x\|_{p,\infty} <\infty\}$ with respect to  $\|\,\|_{p,\infty}$. As an intuitive notation, we will often write 
\[
\set{\md x>\lambda}:=\Ind_{(\lambda,\infty)}(\md x).
\]If $\lambda_{1}+\lambda_{2}=\lambda$,
we have 
\begin{equation}
\tau(\set{\md{x+y}>\lambda})\leq\tau(\set{\md x>\lambda_{1}})+\tau(\set{\md y>\lambda_{2}}).\label{eq:lambda}
\end{equation}
Let $(\Omega, m)$ be a standard $\sigma$-finite measure space. Let $\mathcal{N} := \mathcal{M} \overline{\otimes} L_\infty(\Omega)$ be the von Neumann algebraic tensor product, and let $\varphi := \tau \otimes \int dm$ denote the associated tensor product trace. We will often identify $L_p (\mathcal N)$ with the vector-valued $L_p$-space $L_p (\Omega;L_p (\mathcal{M}))$ for $0<p<\infty$ in the usual way. We will also use the following Hölder type inequalities 
\begin{equation}
\norm{\int_{\Omega}f^{*}g}_{L_r (\mathcal M )}\leq\norm{\p{\int_{\Omega}f^{*}f}^{1/2}}_{L_p (\mathcal M )}\norm{\p{\int_{\Omega}g^{*}g}^{1/2}}_{L_q (\mathcal M )}\label{eq:holder}
\end{equation}
where $0<p,q,r\leq\infty$ be such that $1/r=1/p+1/q$, and $f,g$
are operator-valued functions on $\Omega$ such that the norms at the right hand
side are well-defined and finite (see e.g. \cite[Proposition~1.1]{Mei05}
for more details).

\begin{definition}[Atomic filtration]\label{def:atomicfiltration}
An \textit{atomic filtration}
$ (\mathcal{P}_{n})_{n\geq0}$ on $\Omega$ is a sequence of partitions
of $\Omega$ such that 
\begin{enumerate}
\item[i.] the atoms $A\in \mathcal P _n$ are measurable and have finite measure: 
\[
\forall n\geq0,\ \forall A\in\mathcal{P}_{n},\quad\md A<\infty;
\]

\item[ii.] the partitions are nested: 
\[
\forall m\geq n\geq0,\ \forall A\in\mathcal{P}_{n},\ \exists B\in\mathcal{P}_{m},\quad
A\subset B.
\]
\end{enumerate}
\end{definition}

For an atomic filtration, we can associate a decreasing sequence of von Neumann subalgebras of $\N$, which leads to a natural example of noncommutative martingale. More precisely, we denote by $\N_{k}$ the subalgebra of functions which are constant on atoms of $\mathcal{P}_{k}$, i.e.,
\[
\N_{k}:=\set{\sum_{A\in\mathcal{P}_{k}}f_{A}\Ind_{A} \in \N:f_{A}\in\M}.
\]
We will consider the conditional expectation 
\[
\Eb_{k}:\mathcal N \to \mathcal N _k , \quad f\mapsto\sum_{A\in\mathcal{P}_{k}}\p{\dfrac{1}{\md A}\int_{A}f}\Ind_{A}
\]
that extends to a contraction from $L_{p}(\N)$ onto $L_{p}(\N_{k})$ for
all $p\in[1,\infty]$. Conditional expectations satisfy the bimodule property:
\[ \Eb_{k} (axb) = a\Eb_{k} (x) b,\quad x\in L_{p}(\N),\ a,b\in \N_{k} .\]

Let $E$ designate $L_p$ or $L_{p,\infty}$ and $(f_n)_{n\geq 0}$ a sequence of elements in $L_p(\N)$ for some $p\in (0,\infty]$. Let $(\e_n)_{n\geq 0}$ be a sequence of Rademacher variables on the measure space $[0,1]$. We take the following as definition of noncommutative square function norms
$$
\norm{(f_n)_{n\geq 0}}_{E(\N;\ell_2^{rc})} := \norm{\sum_{n=0}^\infty f_n \otimes \e_n}_{E(\N\overline{\otimes}L_\infty([0,1])},
$$
which is equivalent to the usual square function norm $\norm{(\sum \md{f_n}^2)^{1/2}}_E$ if $\mathcal N$ is commutative. We refer to \cite{Cad18} for more information (where the space $E(\N;\ell_2^{rc})$ was denoted as $\mathscr H _E$). We will also need the notions of noncommutative  (strong and weak) maximal (quasi-)norms denoted by
$$
\norm{(f_n)_{n\geq 0}}_{L_p(\N;\ell_\infty)} = \inf \left\{ \|a\|_{2p} \sup_{n \geq 1} \|y_n\|_\infty \|b\|_{2p}  \mid f_n = ay_n b ,\, a,b\in L_{2p} (\mathcal M ), (y_n ) \subset \mathcal M \right\}
$$and
$$\norm{(f_n)_{n\geq0}}_{\Lambda_{p,\infty}(\N;\ell_\infty)} =\sup_{\lambda>0} \lambda\inf_{\substack{e \in \mathcal{M} \\ \text{projection}}} \left\{ (\tau(e^\perp))^{\frac{1}{p}} : \|ex_ne\|_\infty \leq \lambda \right\}.  
$$
See, for example, \cite{JX07,Hon20ball} for a detailed exposition. These norms weak norms first appeared in Cuculescu's work \cite{Cuc71} to formulate Doob maximal inequality for noncommutative martingales.
\begin{lem}[\cite{Cuc71}]
    We have
    $$\|(\mathbb E _n f)_{n\geq 0} \|_{\Lambda_{1,\infty}(\N;\ell_\infty)} \lesssim \|f\|_1,\quad f\in L_1 (\mathcal{N}).$$
\end{lem}
As can be expected, the noncommutative square function norm dominates the maximal one (see the proof of \cite[Corollary 1.4]{HX21}):
\begin{lem} \label{lem:linftyl2}
For a finite sequence $(f_{n})$
in $L_{1,\infty}(\mathcal{N})$, there exists a universal constant
$c>0$ such that 
$$
\norm{(f_n )}_{\Lambda_{1,\infty}(\N;\ell_\infty)} \leq c\norm{(f_n)}_{L_{1,\infty}(\N;\ell_2^{rc})}.
$$
\end{lem}

\subsection{Amenable groups and asymptotic invariance}

In this section, we recall the characterization of amenable groups through the existence of 
F\o lner sets and discuss two notions of invariance. A more detailed exposition
of similar notions can be found in \cite[Section I]{OW87}.  

Let $G$ be a locally compact second countable group and $r$ a right-invariant
Haar measure on $G$. For a measurable subset $E$ of $G$, we also denote $r(E)$ 
by $\md E$. The most commonly used notion of (almost) invariance is the following. 
Let $K$ be a compact subset of $G$ and $\e > 0$. Assume that $\md{E} < \infty$. 

\begin{definition}
We say that $E$ is \emph{$(\e,K)$-invariant} if 
$$
\md{E  K \setminus E} \leq \e \md{E}\quad \text{and} {\quad \md{E  K^{-1} \setminus E} \leq \e \md{E}}. 
$$
\end{definition}

A related but slightly different notion of invariance can be defined via the consideration
of $K$-boundaries. In the following $E$ will still denote a subset of $G$ of finite
measure and $K$ a compact subset of $G$. The \emph{$K$-boundary} of $E$
denoted by $\partial_{K}(E)$ is defined as the union of the left
translates of $K$: $g  K$ with $g\in G$, such that $g  K$
intersects both $E$ and $E^{c}:=G\setminus E$. The $K$-interior
and $K$-closure of $E$ are defined as 
\[
\mathrm{Int}_{K}(E):=E\setminus\partial_{K}(E)\quad\text{and}\quad\Cl_{K}(E):=E\cup\partial_{K}(E).
\]
\begin{definition}
We will say that $E$ is \emph{$(\e,K)$-boundary-invariant} if 
$$
\md{\partial_K(E)} \leq \e \md{E}. 
$$
\end{definition}
To lighten the notation, we will say that a set is {\it sufficiently invariant} if it is $(\e,K)$-invariant or $(\e,K)$-boundary-invariant for some  large subsets $K$ and sufficiently small constants $\e$ that we do not need to make explicit in the text. 
\begin{rem}\label{rem:inv}	
The notions of $(\e,K)$-invariance and $(\e,K)$-boundary-invariance are nonequivalent but closely related. For convenience of presentation, we keep both notions in this paper. 
The following simple observations can be verified by definition and will be frequently used throughout the paper without further reference:
\begin{enumerate}
    \item if $E$ is $(\e,K \cup \set{e})$-boundary-invariant then $E$ is 
$(\e,K)$-invariant;
\item if $E$ is $(\e,(K^{-1}K)^2)$-invariant then $E  K^{-1}K$ 
is $(\e,K)$-boundary-invariant.
\end{enumerate}
The second verification is slightly less obvious but follows from the identities  
$$
\Cl_K(E) = E K^{-1}K \quad \text{and} \quad E \subset \Int_K(\Cl_K(E)).
$$
To justify the first identity above, note that $ EK^{-1}K $ consists precisely of those translates $ gK $ with $ g \in EK^{-1} $, and in particular $ gK \cap E \neq \emptyset $. Conversely, if a translate $ gK $ intersects $ E $, then there exist $ h \in E $ and $ k \in K $ such that $ g = hk^{-1} $, which implies $ gK \subset EK^{-1}K $. 
For the second inclusion, suppose for contradiction that there exists an element $ h \in gK \cap E $ with $ gK \subset \partial_K(\Cl_K(E)) $. Writing $ h = gk_0 $ for some $ k_0 \in K $, we see that for any $ k \in K $, we have $ gk = hk_0^{-1}k $, which belongs to $ EK^{-1}K $. This contradicts the assumption that $ gK $ lies entirely in the boundary $ \partial_K(\Cl_K(E)) $, since no point of $ EK^{-1}K $ can belong to the boundary. This completes the verification.
\end{rem}


Since we wish to include non-unimodular groups in our construction,
a quick warning is due about left and right invariance. We equipped $G$ 
with a right-invariant measure $r$ and consequently can only expect to have 
almost right-invariant sets. Indeed, by definition of the modular function
$\Delta: G \to \Rb_{>0}$, for any $g\in G$ and $E \subset G$, $\md{gE} = \Delta(g)\md{E}$. 
Hence if $\Delta(g)>1$ and $\e < \Delta(g)-1$, there cannot exist any left-$(\e,\set{g})$-invariant set.  

This notion of invariance allows us to formulate the characterization of amenability that 
we will use in this paper. 

\begin{definition}\label{def:folner}
A locally compact second countable group $G$ is said to be \emph{amenable} if there exists 
a sequence $(F_n)_{n\geq 0}$ of compact subsets of $G$ such that for any $\e >0$ and $K \subset G$ compact, 
$F_n$ is $(\e,K)$-boundary-invariant for sufficiently large $n$. Such a sequence is called
a {\it F\o lner sequence}. 
\end{definition}
Our definition of a Følner sequence differs slightly from conventional formulations found in the literature, which typically employ the standard $(\e , K)$-invariance condition. However, as explained in  Remark \ref{rem:inv}, they yield equivalent characterizations of amenability.  We adopt this version of the definition to avoid unnecessary technical details and notational complications in subsequent analyses. One may always work with the other version of Følner sequences with appropriate minor modifications in later arguments.

\subsection{Reduction to maximal inequalities}\label{subsec:schemaau}

In the following, we present the general recipe for proving Theorem \ref{thm:intromain}. This recipe has previously been used in \cite{HLW21} and originates from the decade-old tradition of studying pointwise ergodic theorems, both for amenable group actions and on noncommutative measure spaces. Most of the ingredients are already gathered in the literature, except for the key maximal inequality discussed in details in later sections. From now on, we consider $\Omega=G$ and $\mathcal N = L_\infty (G) \overline{\otimes}\mathcal M$ for a semifinite von Neumann algebra $\mathcal M$. Recall that $\alpha$ is an action of $G$ on $\mathcal M$ and $(A _n)_{n\geq 0}$ are the ergodic averages associated with a Følner sequence $(F_n)_{n\geq 0}$ in $G$:
$$
A_{n}(x) := \dfrac1{r(F_n)}\int_{F_n} \alpha_g(x) dr(g).
$$

\subsubsection*{\emph{1)} Convergence on a dense subset} By the von Neumann mean ergodic theorem for amenable groups (see e.g. \cite[Theorem 6.7]{Nev06}), it is well known that the subspaces
$$ \{x\in L_2 (\mathcal M ) : \alpha _g x=x,\forall g\in G\}\quad 
 \text{and}\quad
 \mathrm{span} \{ x-\alpha_g x : x\in L_2 (\mathcal{M}) \cap \mathcal M,g\in G\}$$span a dense subspace of $L_2 (\mathcal M)$. It is standard to see from the Følner property of $(F_n)$ that the desired almost
uniform convergence of $(A_n x)$ in Theorem \ref{thm:intromain}  holds for the elements $x$ in $L_1(\M) \cap L_\infty(\M) \subset L_2(\M)$. This observation will be sufficient for later purposes.

\subsubsection*{\emph{2)} Banach principle} 
This principle connects maximal inequalities and almost uniform convergence. To state the principle let us first recall the following notion of almost uniform convergence introduced in \cite{Lan76}, which is well known to be an appropriate analogue of the almost everywhere convergence for the noncommutative setting. \begin{definition}[Almost uniform convergence]\label{def:au} Let $\M$ be a
von Neumann algebra and $\tau$ a normal semifinite faithful trace on $\M$. A sequence of elements
$(x_{n})\subset L_p (\mathcal M )$ is said to converge \textit{bilaterally almost uniformly} (abbreviated as \emph{b.a.u.}) to $x$ if
for any $\e>0$ there exists a projection $e\in\M$ such that $\tau(1-e)\leq\e$
and 
\[
\norm{e(x_{n}-x)e}_{\infty}\to0.
\]
It is said to converge \textit{almost uniformly} (abbreviated as \emph{a.u.}) to $x$ if
for any $\e>0$ there exists a projection $e\in\M$ such that $\tau(1-e)\leq\e$
and 
\[
\norm{(x_{n}-x)e}_{\infty}\to0.
\]
\end{definition}
The following noncommutative Banach principle is well known to experts. A sequence of maps $(\Phi_n)_{n\geq0}$ on $L_p (\mathcal M) $ is of \emph{weak type} $(p,p)\, (p<\infty)$ if $x\mapsto (\Phi_n x)$ defines a bounded map from $L_p (\mathcal M)$ to $\Lambda_{p,\infty} (\mathcal M ;\ell_\infty)$, and is of \emph{(strong) type} $(p,p)$ if it is bounded from $L_p (\mathcal M)$ to $L_{p} (\mathcal M ;\ell_\infty)$. The strong type $(p,p)$ implies the weak one.

\begin{lem}[{see e.g. \cite[Theorem 3.1]{CL21}}] \label{lem:Banachprinciple}Let $1\leq p<2$ (resp. $2\leq p < \infty$) and $(\Phi_{n})_{n\geq0}$ be a sequence of positive linear maps on $L_{p}(\mathcal{M})$ of weak type  $(p,p)$. Then the space of the elements $x\in L_{p}(\mathcal{M})$ such that $\Phi_{n}(x)$ converges b.a.u. (resp. a.u.) is closed in $L_{p}(\mathcal{M})$.
\end{lem}
Together with 1), we see that to prove Theorem \ref{thm:intromain}, it suffices to show that  $(A_n)_{n\geq 0}$ is of type $(p,p)$ for $p>1$ and of weak type $(1,1)$. Indeed, if this is the case, then by the above lemma, the sequence $A_n x$ converges a.u. for all $x\in L_2 (\mathcal M)$; since $L_2 (\mathcal M) \cap L_p (\mathcal M)$ is dense in $L_p (\mathcal M )$, an application of the above lemma again will yield Theorem \ref{thm:intromain}. 
\subsubsection*{\emph{3)} Transference principle.} The classical Calder\'on transference principle states that the maximal inequality for any action can be derived from the maximal inequality for the action of $G$ onto itself by translation. A noncommutative analogue of this principle was recently given in \cite[Theorem~3.1]{HLW21} and \cite[Theorem~3.3]{HLW21}, where the usual translation action on $G$ is replaced by the translation action on the operator-valued function space $\mathcal N = L_\infty (G) \overline{\otimes}\mathcal M$. More precisely, consider the family of operators $(\Ab_n)_{n\geq 0}$ on $L_1(\N) + \N$ defined by 
$$
[\Ab_n(f)](g) = \dfrac1{\md{F_n}}\int_{F_n} f(gh)dm(h),\quad f\in L_1(\N) + \N,\ g \in G,
$$and the principle in  \cite[Theorem 3.1]{HLW21} and \cite[Theorem 3.3]{HLW21} states that the family of maps $( A _n)_{n\geq 0}$  for a general action $\alpha$ is of weak type $(1,1)$ (resp. strong type $(p,p)$) if so is  $(\Ab_n)_{n\geq 0}$. 

\subsubsection*{\emph{4)}  Maximal inequality.} Since the averages $(\Ab_n)_{n\geq 0}$ are contractions on $L_\infty$, the noncommutative Marcinkiewicz-type interpolation theorem \cite[Theorem 3.1]{JX07} yields that to establish the maximal inequalities of type $(p,p)$ for   $(\Ab_n)_{n\geq 0}$, $1\leq p < \infty$, it suffices to prove the one of weak type $(1,1)$. Therefore, together with Observations 1)-3), we see that the proof of Theorem \ref{thm:intromain} reduces to the study of the following result, which is our main objective. 
\begin{thm}\label{thm:maxmain}
Let $G$ be a locally compact second countable amenable group. There exists a F\o lner sequence $(F_n)_{n\geq 0}$ such that   $(\Ab_n)_{n\geq 0}$ is of weak type $(1,1)$, i.e., for any $f \in L_{1}(\N)$,
$$
\norm{(\Ab_n(f))_{n\geq0}}_{\Lambda_{1,\infty}(\N;\ell_\infty)} \les \norm{f}_1.
$$
\end{thm}
 
The remainder of this paper will be devoted to the proof of Theorem \ref{thm:maxmain}.
As mentioned previously,  the methods used in \cite{HLW21} do not apply to the actions of general amenable groups without doubling conditions. Note that Lindenstrauss's proof of the commutative version of Theorem \ref{thm:maxmain} (see \cite{Lin01}) does not extend directly to the noncommutative setting neither due to the current lack of understanding of noncommutative covering lemmas. We will explore another path based on noncommutative variational inequalities and Ornstein-Weiss tilings, as explained below. 

\subsection{Outline of the strategy and main results}

As mentioned previously, our approach will be based on the existence of F\o lner sequences equipped with certain {\it admissible} filtrations. We take inspiration from the Euclidean setting where $G = \Rb^d$, $\N =L_\infty (\mathbb R ^d )\overline{\otimes} \M $, and the following are considered:
\begin{itemize}
    \item balls with dyadic radii $F_n \coloneqq B(0,2^n)$,  $n\geq 0$;
    \item dyadic filtration  $\Pc_n := \set{[0,2^n)^d + 2^nk : k\in \Zb^d}$, $n\geq 0$;
    \item dyadic cubes  $B_n := [0,2^n]^d$.
\end{itemize}
The strategy developed in \cite{JRW03} and extended to operator-valued functions in \cite{HX21,Xu21} to obtain information on the ergodic averages is to compare them with the conditional expectations coming from the filtration $(\Pc_n)_{n\geq 0}$. A posteriori, we observe that the key geometric properties implicitly used in these works are the following:
\begin{itemize}
\item for any $n\geq k$, $F_n$ is $(c2^{k-n},B_k)$-boundary-invariant, 
\item for any $k \geq n$ and $A \in \Pc_k$, $A$ is $(c2^{n-k},F_k)$-invariant, 
\end{itemize}
where $c>0$ is a harmless constant. Returning to a general locally compact second countable amenable group $G$, we would like to build similar structures on an intuitively large portion of $G$.  To quantify how large a subset covers $G$, we are going to use the notion of \emph{(lower) Banach density} introduced in \cite{BBF10,DHZ19}: for $ E\subset G$, we define
\[ 
\underline{\mathrm B} (E) = \sup_{K \subset G} \inf_{g\in G} \dfrac{\md{g  K \cap E}}{\md{gK}},
\]
where the supremum is taken over all compact subsets $K$ of $G$. Taking all of the above into consideration, we are led to study the
following natural notions.

\begin{definition}\label{def:filtered}
A \emph{filtered} F\o lner sequence is a triple $\FS = ((F_n)_{n\geq 0},(B_n)_{n\geq 0},(\mathcal{P}_n)_{n\geq 0})$ where 
$(F_n)_{n\geq 0}$ is a F\o lner sequence,  $(B_n)_{n\geq 0}$ are sequences of compact subsets of $G$, and $(\mathcal{P}_n)_{n\geq 0}$ is an atomic filtration such that 
\begin{enumerate}
\item[i.] for any $n\geq 0$ and $A\in \Pc_n$, there exists $\gamma\in G$ such that $A \subset \gamma  B_n$,
\item[ii.] for any $n\geq k$, $F_n$ is $(2^{k-n},B_k)$-boundary-invariant.
\end{enumerate}
For such a sequence, we say that an atom $A \in \Pc_k$ is \emph{admissible} if  $A$ is $(2^{n-k},F_n)$-invariant for any $ n<k$. Denote by $\Pc_k^a$  the family of all admissible atoms of $\Pc_k$ and write $\Adm(\S):=\bigcap_{k\geq0}\bigcup_{A\in\Pc_{k}^{a}}A$. We say that $\FS$ is \emph{$c$-admissible} for a constant $c\in (0,1)$ (or simply \emph{admissible}) if $\underline{\mathrm{B}}(\Adm(\S) ) \geq 1-c$. A function defined on $G$ is said to be \emph{admissible} if it is supported in $\Adm(\S)$.
\end{definition}

By a slight abuse of terminology, we say that a F\o lner sequence $(F_n)_{n\geq 0}$ is an {\it ($c$)-admissible F\o lner sequence} as well if it can be completed into a triple $((F_n)_{n\geq 0},(B_n)_{n\geq 0},(\Pc_n)_{n\geq 0})$ satisfying the conditions above. Recall the ergodic averages $(\Ab_n)_{n\geq 0}$ associated with $(F_n)_{n\geq 0}$ as well as the conditional expectations $(\Eb_n)_{n\geq 0}$ associated with $(\Pc_n)_{n\geq 0}$. As in the aforementioned work  \cite{JRW03, HX21, Xu21}, we will prove a square function estimate for the associated difference operator 
\begin{equation}\label{eq:diff}
    \Db := (\Ab_n - \Eb_n)_{n\geq 0}.
\end{equation}
Two questions arise:

(a) Can admissible sequences be constructed beyond the Euclidean case? 

(b) Are the imposed conditions sufficient to bound the square function $\Db$?
  
Concerning Question (a), various natural and explicit examples can be found: in the context of Heisenberg groups, lamplighter groups, and Brieussel-Zheng's diagonal products (see e.g. \cite{DKLT22,Esc24}). The lamplighter group is a typical example that does not admit doubling F\o lner sequences (see \cite[Corollary 5.5]{Lin01}), it eludes the techniques of \cite{HLW21}. In general, the existence of filtered F\o lner sequences is strongly related to the {\it tiling} properties of amenable groups, which already had a deep impact in ergodic theory. Inspired by \cite{OW87} and \cite{DHZ19}, we will prove the following. 

\begin{mainthm}\label{thm:wellfiltered}
Let $G$ be a locally compact second countable amenable group. Any F\o lner sequence in $G$ admits
an admissible F\o lner subsequence. 
\end{mainthm}

To answer Question (b), we follow the strategy of \cite{JRW03, HX21, Xu21} combined with the non-doubling noncommutative Calder\'on-Zygmund decomposition defined in \cite{CCP21}.

\begin{mainthm}\label{thm:maindiffop}
Let $\Db$ be the difference operator defined in \eqref{eq:diff} associated with an admissible F\o lner sequence $(F_n )_{n\geq0}$. Then for any $f\in L_1(\N)$,
$$
\norm{\Db f}_{L_{1,\infty}(\N;\ell_2^{rc})} \les \norm{f}_1.
$$
\end{mainthm}

Theorem \ref{thm:wellfiltered} is proved in Section \ref{sec:geometry} and Theorem \ref{thm:maindiffop} is proved in Sections \ref{sec:CZdec} and \ref{sec:diffop}. Combining the two gives us a proof of Theorem \ref{thm:maxmain}. 

\begin{proof}[Proof of Theorem \ref{thm:maxmain}]
By Theorem \ref{thm:wellfiltered}, we know that $G$ admits an admissible F\o lner sequence $((F_n)_{n\geq 0},(B_n)_{n\geq 0}, 
(\Pc_n)_{n\geq 0})$. Let $f \in L_1(\N)$. By the quasi-triangle inequality, 
$$
\norm{(\Ab_n f)_{n\geq 0}}_{\Lambda_{1,\infty}(\N;\ell_\infty)}
\les \norm{(\Eb_n f)_{n\geq 0}}_{\Lambda_{1,\infty}(\N;\ell_\infty)}
+ \norm{(\Eb_n f - \Ab_n f)_{n\geq 0}}_{\Lambda_{1,\infty}(\N;\ell_\infty)}.
$$
By Cuculescu's inequality \cite{Cuc71}, the first term is dominated by $\norm{f}_1$. For the second term, we can use the fact that maximal functions are dominated by square functions (Lemma \ref{lem:linftyl2}) and, since $(F_n)_{n\geq 0}$ is admissible, apply Theorem \ref{thm:maindiffop} to obtain:
\[
\norm{(\Eb_n f - \Ab_n f)_{n\geq 0}}_{\Lambda_{1,\infty}(\N;\ell_\infty)}
\les \norm{(\Eb_n f- \Ab_n f)_{n\geq 0}}_{L_{1,\infty}(\N;\ell_2^{rc})}
\les \norm{f}_1.  \qedhere
 \] 
\end{proof}
As a result of Theorem \ref{thm:maxmain} and the framework outlined in Subsection~\ref{subsec:schemaau}, we can conclude with the following theorem, which gives a precise answer to the original motivating question.
\begin{thm}\label{thm:main}
        Let $G$ be a locally compact second
countable amenable group and $(F_{n})$ be an admissible Følner
sequence on $G$ (which always exists according to Theorem \ref{thm:wellfiltered}). If $\alpha$ is a w{*}-continuous action of $G$ on a semifinite von
Neumann algebra $(\mathcal{M},\tau)$ by $\tau$-preserving automorphisms and we consider the averages
$$ A _n x=\frac{1}{|F_{n}|}\int_{F_{n}}\alpha_{g}xdm(g),\quad x\in \mathcal M,$$
then  $( A _n )$ extends to a sequence of maps of weak type $(1,1)$ and strong type $(p,p)$ ($1<p\leq\infty$). For all $x\in L_{p}(\mathcal{M})$ with $1\leq p<2$ (resp. $2\leq p < \infty$), 
\[
 A _n x \to Px\ \text{b.a.u. (resp. a.u.)},
\]
where $P$ is the projection onto the space of fixed points
of $\alpha$. 
\end{thm}
This theorem in particular applies to the aforementioned concrete examples of canonical F\o lner sequences with tilings on Heisenberg groups, lamplighter groups, and Brieussel-Zheng's diagonal products (\cite[Propositions 6.15 and 6.19]{DKLT22}, \cite[Theorem 3.9]{Esc24}).
\subsection{Further remarks}
\label{sub:furtherrem}

In the classical setting, our method provides new tools for the study of quantitative pointwise ergodic theorems for amenable groups. As an illustration, let us briefly present a new pointwise ergodic theorem on $L_1$. Given a sequence of measurable functions $(f_n)_{n\geq 0}$ on a measure space $X$ which converges almost everywhere, it is natural to consider the local quantitative behavior of the Cauchy sequence $(f_n (x))_{n\geq 0}$  for a.e. $x\in X $. More precisely, for any $\lambda>0$ and  $x\in X $, the \emph{jump quantity} $N_\lambda((f_n)_{n\geq 0})(x)$ is defined as the maximal number $N$ for which there are $n_1 <n_2<\cdots<n_N$ such that
\[ |f_{n_i}(x) - f_{n_{i+1}}(x)| >\lambda,\quad 1\leq i \leq N-1.  \]

\begin{thm}\label{thm:quantitative}
Let $G$ be a locally compact second
countable amenable group and $(F_{n})$ be an admissible Følner
sequence on $G$ (which always exists according to Theorem \ref{thm:wellfiltered}). Consider a measure-preserving action of $G$ on a measure space $(X,\mu)$ and let 
\[  A _n f =\frac{1}{\md{F_{n}}} \int_{F_{n}}s. f\ dm(s),\quad f\in L_1 (X). \]
Then there is a constant $c>0$ such that for all $\lambda>0$ and $m>0$,
\[  \mu \{  {N_\lambda((  A _n f)_{n\geq 0})}>m\}\leq \frac{c}{\lambda\sqrt{m}}\|f\|_1 ,\quad f\in L_1 (X). \]  
\end{thm}
\begin{proof}
By \cite[Theorem 1.7]{HL21quantitative}, it suffices to consider the case $X=G$ with the translation action of $G$ on itself (in \cite{HL21quantitative} only the ball averages with respect to an invariant metric are considered, but the proof applies to general Følner averages as well). As in \cite{HL21quantitative}, we note that
\begin{align*}
	\sqrt{N_\lambda((\mathbb A _n f)_{n\geq 0})}&\leq   \sqrt{N_{\lambda/2}((\mathbb D _n f)_{n\geq 0})} + \sqrt{N_{\lambda/2}((\mathbb E _n f)_{n\geq 0})}\\
	&	\leq \dfrac{ {2}}{  \lambda} \left(\sum_n |\mathbb D_n f |^2\right)^{1/2} + \sqrt{N_{\lambda/2}((\mathbb E _n f)_{n\geq 0})}.
\end{align*} 
The weak (1,1) inequality for the first term follows from Theorem \ref{thm:maindiffop}, and the one for the second is established in \cite[Theorem 6.8]{JKRM88osc} (and implicitly in \cite{PX88}).
\end{proof}

Before ending the section, we remark that our method is also useful for the study of some other more general or stronger form of ergodic theorems, for instance, the square estimates $(\mathbb{A}_n -\mathbb{A}_{n-1})_{n\geq 0}$ and differential transforms similar to \cite{HX21,Xu21}, maximal and individual ergodic theorems for group actions on noncommutative $L_p$-spaces with a fixed $p\in (1,\infty)$ as in \cite{HLW21}. Our previous results for noncommutative $L_p$-spaces with $p\in (1,\infty)$ hold
true as well for general non-tracial von Neumann algebras, up to standard adaptations as in \cite{JX07,haagerupjx10reduction}.




\section{Existence of admissible Følner sequences}
\label{sec:geometry}

This section is dedicated to proving Theorem \ref{thm:wellfiltered}. The argument will be done by inductive constructions of appropriate quasi-tilings of groups. From now on, we will always assume that $G$ is a locally compact second countable amenable group. A {\it quasi-tiling} is formally a collection of subsets $(T_i)_{i\in I}$ that we are going to endow with some useful structure.  

\begin{definition}[Quasi-tiling]
    A quasi-tiling is a triple $\T = (I,(T_i)_{i\in I},(\gamma_i)_{i\in I})$ where $I$ is a well-ordered set, $(T_i)_{i\in I}$ is a family of subsets of $G$ and $(\gamma_i)_{i\in I}$ is a family of elements of $G$. The element $\gamma_i^{-1}$ plays the role of the center of the tile $T_i$. 

    We say that a quasi-tiling is:
\begin{itemize}
    \item {\bf bounded} if $\bigcup_{i\in I} \gamma_i  T_i$ is precompact,
    \item {\bf based on $\S$} for a family $\mathcal S$ of subsets of $G$ if $\set{\gamma_i  T_i}_{i\in I} \subset \S$, 
    \item {\bf $\e$-disjoint} if for any $i\in I$, 
    $$
    \md{ T_i \cap \left(\bigcup_{j < i} T_j \right)} 
< \e \md{ T_i}.
    $$
    \item {\bf disjoint} if the $T_i$'s are disjoint,
    \item {\bf $(\delta,K)$-(boundary)-invariant} if each $T_i$ is $(\delta,K)$-(boundary)-invariant.
\end{itemize}
We will write $\underline{\mathrm B}(\T) := \underline{\mathrm B}(\bigcup_{i\in I} T_i)$ and more generally by abuse of notation $\T := \bigcup_{i\in I} T_i$. 
\end{definition}

Our aim is to construct quasi-tilings of $G$ that maintain high Banach density while being disjoint, bounded, and sufficiently boundary-invariant. For discrete groups, such tiling properties were established by Downarowicz, Huczek, and Zhang in \cite[Theorem 5.2]{DHZ19}, where the constructed tilings are even exact (i.e. $G=\bigcup_{i\in I} T_i$) and further satisfy a zero-entropy condition. The foundational framework for these results traces back to the seminal work of Ornstein and Weiss \cite{OW87}. Our approach synthesizes techniques from \cite{OW87} and \cite{DHZ19} with a generalization of the core argument to non-unimodular groups. In the appendix we will also establish the exact tiling property for unimodular locally compact amenable groups, with a non-trivial extension of the marriage lemma methodology of \cite{DHZ19}  to the non-discrete setting.

\subsection{Banach density}
We first study a few preliminary properties of Banach density. Recall that for $ E\subset G$, we define
\[ 
\underline{\mathrm B} (E) = \sup_{K \subset G} \inf_{g\in G} \dfrac{\md{g  K \cap E}}{\md{gK}},
\]
where the supremum is taken over all compact subsets $K$ of $G$. In the following we prove that this supremum can be equivalently attained by considering  the compact subsets $K$ that are sufficiently invariant. As far as we know, this equivalent formulation was previously shown only for discrete groups  \cite{BBF10,DHZ19} and the following is a generalization to the non-discrete setting.

\begin{prop}\label{prop:banachdensity}
Let  $(F_n)_{n\geq 0}$ be a F\o lner sequence on $G$. We have
$$
\bd{E} = \limsup_{n\to\infty} \inf_{g\in G} \dfrac{\md{g  F_n \cap E}}{\md{g  F_n}}.
$$
\end{prop}

To prove this proposition, let us isolate the following key computation.

\begin{lem}\label{lem:ineqBdensity}
    Let $K,A$ and $D$ be precompact measurable subsets of $G$. Then
 \begin{equation}\label{eq:ineqBdensity}
      \md{K}\md{D}\inf_{g\in D} \dfrac{\md{A \cap g  K}}{\md{g  K}} \leq 
\int_D \md{A \cap g  K}dl(g) 
\leq \md{K}\md{A\cap D  K}.
 \end{equation}

\end{lem}

\begin{proof}
The first inequality follows from the relation between the Haar measures and modular functions. More precisely, we have
\begin{align*}
\int_D \md{A \cap g  K}dl(g) &\geq \inf_{g\in D} \dfrac{\md{A \cap g  K}}{\md{g  K}} \int_D \md{g  K} dl(g) \\
&= \inf_{g\in D} \dfrac{\md{A \cap g  K}}{\md{g  K}} \md{K} \int_D \Delta(g)dl(g) \\
&= \md{K}\md{D}\inf_{g\in D} \dfrac{\md{A \cap g  K}}{\md{g  K}}.
\end{align*}
For the second inequality, notice first that for $g\in D$, we have $g  K \subset  D  K$. Hence, 
$$
\int_D \md{g  K \cap A} dl(g) = \int_D \md{g  K \cap A \cap D  K}dl(g) 
\leq \int_G \md{g  K \cap A \cap D  K}dl(g).
$$
Setting $E = A \cap D  K$ and applying Fubini's theorem, we obtain
\begin{align*}
    \int_G \md{g  K \cap E}dl(g)
    &= \int_{g\in G} \int_{h\in E} \Ind_{g  K}(h) dr(h) dl(g) \\
    &= \int_{h\in E} \int_{g\in G} \Ind_{h  K^{-1}}(g) dl(g) dr(h) \\
    &= \md{E}l(K^{-1}) = \md{E}\md{K}.  \qedhere
\end{align*}
\end{proof}

\begin{cor}\label{cor:tech:goodtranslate}
    Let $K,A$ and $D$ be compact subsets of $G$.  Assume that $D$ is $(\e,K)$-invariant. Then 
\[
\inf_{g\in D} \dfrac{\md{A \cap g  K}}{\md{g  K}} \leq \dfrac{\md{A \cap D}}{\md{D}} + \e.
\]
\end{cor}

\begin{proof}
Note that if $D$ is $(\e,K)$-invariant then 
$$
\md{A\cap D  K} \leq \md{A \cap D} + \md{D  K \setminus D} \leq \md{A \cap D} + \e \md{D}.
$$
By substituting the right-hand side of \eqref{eq:ineqBdensity} with the above expression, we obtain the desired inequality. 
\end{proof}

\begin{proof}[Proof of Proposition \ref{prop:banachdensity}]
The inequality
\[
\bd{E} \geq \limsup_{n\to\infty} \inf_{g\in G} \frac{\md{g  F_n \cap E}}{\md{g  F_n}}
\]
is clear from the definition of $\bd{E}$. To establish the reverse inequality, consider a compact subset $K \subset G$. Let $\e  > 0$. Note that if a set $F$ is $(\e , K)$-invariant, then $g  F$ remains $(\e , K)$-invariant for any $g \in G$. Consequently, by Corollary \ref{cor:tech:goodtranslate}, we obtain
\[
\inf_{g\in G} \frac{\md{g  K \cap E}}{\md{g  K}}
\leq \inf_{g'\in G}\inf_{g\in g'F} \frac{\md{g  K \cap E}}{\md{g  K}}
\leq \inf_{g'\in G}\frac{\md{g'  F \cap E}}{\md{g'   F}} + \e.
\]
Applying this inequality to a Følner sequence $(F_n)_{n\geq 0}$ with $n$ tending to infinity proves the Proposition. 
\end{proof}
\begin{rem}
    We see from the above proof that it suffices to consider in the proposition a sequence $(F_n)$ of subsets which are eventually $(\e ,K)$-invariant for all $\e$ and $K$ instead of being  $(\e ,K)$-boundary-invariant in Definition \ref{def:folner}.
\end{rem}
We conclude this subsection with a simple observation on the Banach density of intersections of subsets.
\begin{lem}\label{lem:banachunion}
    Let $E_1$ and $E_2$ be two subsets of $G$. Then \[ \underline{\mathrm{ B}}(E_1 \cap E_2  )\geq \underline{\mathrm{ B}}(E_1   )+\underline{\mathrm{ B}}(  E_2  )-1.\]
\end{lem}
\begin{proof}
    The proof is a straightforward verification by definition. Note that a subsequence of a F\o lner sequence is still a F\o lner sequence. By passing to subsequences several times, we may find a distinguished F\o lner sequence $(F_n)$ such that the superior limit in Proposition \ref{prop:banachdensity}
 is attained by  $(F_n)$  for $E_1 ,E_2$ and $E_1 \cap E_2$ simutaneously, i.e.,    \[ \bd{E_1} = \lim_{n\to\infty} \inf_{g\in G} \dfrac{\md{g F_n \cap E_1}}{\md{g  F_n}}, \quad \bd{E_2} = \lim_{n\to\infty} \inf_{g\in G} \dfrac{\md{g  F_n \cap E_2}}{\md{g F_n}} \]and \[ \bd{E_1\cap E_2} = \lim_{n\to\infty} \inf_{g\in G} \dfrac{\md{g  F_n \cap E_1\cap E_2}}{\md{g  F_n}}.\]Note that for any compact subset $K\subset G$, we have \[|K\cap E_1 \cap E_2 | =  |(K\cap E_1 ) \setminus (K\setminus E_2 ) |\geq |K\cap E_1 | - (|K| - |K\cap  E_2|).\]So the desired inequality follows.
\end{proof}

\subsection{Maximal quasi-tilings} Following \cite{OW87, DHZ19}, the first step of the construction is to tile $G$ with translates of a finite set of shapes $\S$. A way to formalize this idea is to consider maximal families of translates that $\e$-intersect.

\begin{definition}[Order and union] Les $\T^1 = (I^1,T^1,\gamma^1)$ and $\T^2 = (I^2,T^2,\gamma^2)$ be two quasi-tilings.

Define the following order: $\T^1 \succ \T^2$ if $I^2$ is an initial segment of $I^1$ (as ordinals) and the canonical increasing map $f:I^2 \to I^1$ satisfies, for any $i\in I^2$, $T^2_i = T^1_{f(i)}$ and $\gamma^2_i = \gamma^1_{f(i)}$.

The union of two tilings $\T^1 \sqcup \T^2$ is defined by concatenation (meaning in particular that $\T^1 \sqcup \T^2 \neq \T^2 \sqcup \T^1$ since their elements are ordered differently). 

The limit of a non-decreasing family of tilings $(\T^j)_{j\in J}$ is defined as $\T = (I,T,\gamma)$ where $I = \sup_{j\in J} I^j$ (as ordinals). Let $f_j$ be the canonical map from $I_j$ to $I$. For any $i\in I$, there exists $j\in J$ and $i_j \in I^j$ such that $i = f_j(i_j)$ and we define $T_i = T_{i_j}$. This definition does not depend on the choice of $j$ by definition of the order on quasi-tilings. 
\end{definition}

Note that the properties of being $B$-bounded, $\S$-based, or $\e$-disjoint are preserved by increasing limits so that by Zorn's lemma we can consider maximal quasi-tilings satisfying those properties. This allows to formulate the following key observation.   As mentioned previously, by abuse of notation we write $\T := \bigcup_{i\in I} T_i$ and $\R := \bigcup_{j\in J} R_j$.

\begin{lem}\label{lem:maximal1}
Let $K$ be a compact subset of $G$. Let $\T$ be a maximal, $\e$-disjoint, $\set{K}$-based quasi-tiling of $G$. Then $\bd{\T} \geq \e$.
\end{lem}

\begin{proof}
By maximality, for any $g\in G$, $\md{gK \cap \T} \geq \e\md{gK}$. Hence, by definition of the lower Banach densito, $\bd{\T} \geq \e$. 
\end{proof}

With this lemma alone we cannot hope to obtain good covering properties with tiles that are almost disjoint. To go further, the idea first exploited in \cite{OW87} is then to cover $G\setminus \T$ with smaller tiles to get a higher Banach density. 

\begin{lem}\label{lem:maximal2}
Let $0<\e < 1/4$. Let $K$ be a compact subset of $G$ and $\R = (J,R,\alpha)$ be an {$(\e^2 ,K)$}-boundary invariant, $\e$-disjoint, $B$-bounded quasi-tiling of $G$. Let $\T$ be a maximal, $\e$-disjoint, $\set{K}$-based quasi-tiling disjoint from $\R$ ($\R \cap \T = \emptyset$). Then 
$$
1  - \bd{\R \sqcup \T} \leq \max\p{\e, (1-\e/3)(1-\bd{\R})}.
$$
\end{lem}

\begin{proof}
    We first show that if $D$ is a sufficiently invariant compact subset of $G$, then
$$
    |D\setminus (\mathcal R \sqcup \mathcal T ) | \leq \max\p{\e\md{D}, (1-\e/3)\md{D\setminus\R}}.
$$    
    More precisely, we will show that this inequality is valid as long as $D$ is $(\e^2,K)$-boundary-invariant and $(\e,B)$-boundary-invariant.
    
    Set $D' = D\setminus \R$. Note that if $\md{D'}   \leq \e\md{D}$, then $|D\setminus (\mathcal R \sqcup \mathcal T )|/|D| \leq \varepsilon$ and there is nothing to prove. So assume that $\md{D'} \geq \e\md{D}$. 
    
    By maximality of $\T$, for any $g\in \Int_K(G\setminus \mathcal R)$, 
    $$
    \md{gK \cap \T} \geq \e\md{gK}.
    $$
    Applying Lemma \ref{lem:ineqBdensity} to $\Int_K(D')$, $K$ and $\T$,
    $
    \md{\Int_K(D')}\e \leq \md{\T \cap \Int_K(D')K}.
    $
    We can assume without loss of generality that $K$ contains $e$ so $\Int_K(D')K \subset D'$ and 
      \begin{equation}\label{eq:lem:maximal2}
    \md{\Int_K(D')}\e \leq \md{\T \cap D'}.
    \end{equation}
    Let $J' \subset J$ be the set of indices $j$ such that $R_j \cap D \neq \emptyset$. The $\e$-disjointness of $R_i$'s yields by induction that 
    \[
    \left |\bigcup_{j\in J'} { R}_j \right| \geq    (1-\varepsilon)\sum_{j\in J'}|{ R}_j|.
    \]
    Combining the $(\e^2,K)$-boundary-invariance of $R_j$'s, we have
$$
\md{\partial_{K}(\bigcup_{i\in I}{ R}_i)} 
\leq \sum_{j\in J'} \md{\partial_{K}{ R}_j} 
\leq \e^2\sum_{j\in J'}\md{{ R}_j}
\leq \frac{\e^2}{1-\e }\md{\bigcup_{i\in I}{ R}_j}.
$$
As $\bigcup_{j\in J'} { R}_j \subset \mathrm{Cl}_{B}(D)$ and $D$ is $(\e,B)$-boundary-invariant
$$
\md{\bigcup_{i\in I}{ R}_i} \leq (1+\e)\md{D}.
$$
Combining these estimates we obtain
\begin{align*}
\md{\partial_{K}(D')} 
&\leq \md{\partial_{K}(D)} + \md{\partial_{K}(\bigcup_{j\in J'}{ R}_j )} \\
&\leq \md{D}\e^2(1 +  (1-\e)^{-1}(1+\e)) \\
&\leq \md{D'}\e(1 +  (1-\e)^{-1}(1+\e)) \\
&\leq \frac{8\e }{3}\md{D'}
\end{align*}
where we used on the last line that $\e \leq 1/4$. Hence $\md{\Int_K(D')} \geq (1-8\e/3)\md{D'} \geq \md{D'}/3$ and by \eqref{eq:lem:maximal2}, 
$$
\dfrac{\e\md{D'}}{3} \leq \md{\T \cap \D'}.
$$

So $\md{\D \setminus (\R \cup \T)} = \md{D'} - \md{\T \cap \D'} \leq (1-\e/{ 3})\md{D \setminus \R}$. Hence we have shown that for any $D$ sufficiently invariant:
$$
\md{D\setminus (\R \cup \T)} \leq \max \p{\e\md{\D}, (1-\e/{ 3})\md{\D\setminus \R}}.
$$
We can apply this inequality for $D = gF_n$ for $(F_n)$ a F\o lner sequence, $n$ large enough, and $g$ any element of $G$. Then we conclude by letting $n$ go to infinity and applying Proposition \ref{prop:banachdensity}.
\end{proof}

\subsection{Quasi-tilings} \label{sec:quasitilings} 
This subsection develops the core construction of quasi-tilings. We first analyze the $\e $-disjoint quasi-tilings and then refine the result for exactly disjoint ones.

\begin{prop} \label{prop:fromOW} Let $0<\e < 1/4$ and
  $\S = (S_{1},\dots,S_{N})$ be a family of compact subsets of $G$ containing the neutral element $e$. Assume that 
  \begin{itemize}
  \item $(1-\e/3)^{N}\leq\e$, 
  \item for any $j>i$, $S_{i}$ is $(\e^3,S_{j})$-boundary-invariant. 
  \end{itemize}
  Then there exists an  $ \e $-disjoint quasi-tiling $\T$ based on $\S$ with $\underline{\mathrm B}(\T)\geq 1- \e$. 
\end{prop}

\begin{proof}
    We are going to construct $\T$ by induction as a union $\T = \T^1 \sqcup \dots \sqcup \T^N$ where $\T^i$ is $\set{S_i}$ based. More precisely, define inductively $\T^i$ to be a maximal, $\e$-disjoint, $\set{S_i}$-based quasi-tiling with centers in $\Int_{S_i}(G\setminus \sqcup_{j<i} \T^j).$ Lemma \ref{lem:maximal2} guarantees that for any $2\leq i\leq N$, 
    $$
    1 - \bd{\sqcup_{n\leq i} \T^n} \leq \max\p{\e,(1-\e/3)(1-\bd{\sqcup_{n<i}\T^n}}.
    $$
Hence, by induction and Lemma \ref{lem:maximal1}, 
\[
1 - \bd{\sqcup_{n\leq N} \T^n} \leq \max\p{\e,(1-\e/3)^N} \leq \e.  \qedhere
\]
\end{proof}

\begin{prop}\label{prop:disbound}
Let $\delta >0$ and $K \subset G$ compact. There exists a tiling $\T$ of $G$ with $\bd{\T} > 1-\delta$ and that is disjoint, bounded, and $(\delta,K)$-boundary-invariant.
\end{prop}

\begin{proof}
           Let $\e>0$. Recall that $ \mathrm{Cl}_{K} (E) = EK^{-1} K$ for any compact subset $K\subset G$.  By Proposition~\ref{prop:fromOW} and amenability, we may choose $(S_1,\dots,S_N)$ such that there exists an $\e$-disjoint quasi-tiling $\Tcl$ based on $(\Cl_K(S_1),\dots,\Cl_K(S_N))$ with $\underline{\mathrm B} (\T ')>1-\e$. Assume additionally that the $S_i$'s are $(\e,K)$-boundary-invariant. From the proof of Proposition~\ref{prop:fromOW}, $\Tcl$ can be decomposed as $\Tcl = \sqcup_{n\leq N} \Tcl^n$ where $ \Tcl^n$ is $\set{\Cl_K(S_n)}$-based.

    Write $\Tcl^n = (I^n, \Tl^n, \gamma^n)$ and define $(\T')^n = (I^n, T'^n, \gamma^n)$ with 
    $$
    T'^n_i := \gamma^n_iS_n.
    $$
Set $\T' = \sqcup_{n\leq N} (\T')^n =: (I, T', \gamma)$ and finally $\T := (I, T, \gamma)$ with 
$$
T_i := T'_i \setminus \p{\bigcup_{j<i} T'_j}.
$$
Note that $\T$ is disjoint by construction and bounded since for any $n\leq N$ and $i\in I_n$,
    $$
    \gamma_i^{-1}T_i \subset \bigcup_{m\leq N} S_m
  .  $$  
   Let $i \in I$. We will show that $T_i$ is $( 3\e , K)$-boundary-invariant. Since by definition $T_i = T_i' \setminus \p{\bigcup_{j<i} T_j'}$, we have
    $$
    \partial_K(T_i) \subset \partial_K(T_i') \cup \p{\bigcup_{j<i} \partial_K(T_j') \cap \Cl_K(T_i)}.
    $$
    We estimate both terms separately. For the first term, by the $(\e,K)$-invariance of the $S_n$'s (and consequently of $T_i'$), we have 
    $$
    \md{\partial_K(T_i')} \leq \e \md{T_i'}.
    $$
Note that by the 
$\e$-disjointedness of 
$\T'$, we have \[\md{\Cl_K(T_i') \cap \p{\bigcup_{j<i}\Cl_K (T_j')}} \leq \e \md{\Cl_K(T_i')} \leq \e (1 + \e) \md{T_i'}.\]So  \begin{align*}
    \md{T_i'} - \md{T_i} &= \md{T_i' \cap \p{\bigcup_{j<i}T_j'}}\\
    &\leq \md{\Cl_K(T_i') \cap \p{\bigcup_{j<i}\Cl_K (T_j')}}\\
    &\leq \e (1 + \e) \md{T_i'} 
    \end{align*}
and hence \begin{equation} \label{eq:ti}
    \md{T_i} \geq (1-2\e)\md{T_i'}.
\end{equation} 
Moreover, for the second term we have
    $$
    \md{\p{\bigcup_{j<i} \partial_K(T_j') \cap \Cl_K(T_i)}} 
    \leq \md{\p{\bigcup_{j<i} \Cl_K(T_j') \cap \Cl_K(T_i')}}
    \leq \e (1+\e) \md{T_i'}. 
    $$
Combing the two estimates yields the desired  $( 3\e , K)$-boundary-invariance of $T_i$.   

We now turn to the Banach density of $\T$. Set $B = \bigcup_{n\leq N} \Cl_K(S_n)$ and let $D$ be an $(\e,B)$-boundary-invariant compact subset of $G$. Using Ineq. \eqref{eq:ti} above, we see that   $$
     \md{D \cap \bigcup_{i \in I}T_i }
    \geq \sum_{T_i \subset D} \md{T_i}   
    \geq (1-2\e)\sum_{ T_i \subset D} \md{ T_i '}.
    $$Note that \[\sum_{ T_i \subset D} \md{ T_i '} \geq \frac{1}{1+\e}\sum_{ T_i \subset D}   |\mathrm{Cl}_K (T_i ') | \geq \frac{1}{1+\e}  |\bigcup_{ T_i \subset D} \mathrm{Cl}_K (T_i ') | \geq 
    \frac{1}{1+\e}  \left|\left(D\cap \bigcup_{ i\in I} \mathrm{Cl}_K (T_i ') \right) \setminus \partial_B (D) \right|. \]
So, using the Banach density $\bd{\T'}$ and the  $(\e,B)$-boundary-invariance of $D$, we obtain that \[ \md{D \cap \bigcup_{i \in I}T_i } \geq   \frac{1-2\e}{1+\e} |D| ,  \]
which implies $\underline{\mathrm B} (\T) \geq  \frac{1-2\e}{1+\e} $. Note that $\e$ was arbitrarily chosen. Thus the proof is complete.
\end{proof}
\subsection{Construction of atomic filtration}
Using the construction from the previous subsection, we can now create the atomic filtration which produces an admissible F\o lner subsequence. This will conclude the proof of Theorem \ref{thm:wellfiltered}.

\begin{lem}\label{lem:tilingforfiltration}
Let $\R = (I,R,\alpha)$ be a $B$-bounded tiling of $G$ with $\underline{\mathrm{B}}(\R) \geq 1-c$ for some $c>0$. Let $\delta>0$ and $K\subset G$ be compact. Then there exists a $(\delta,K)$-invariant disjoint bounded tiling $\T = (J,T,\gamma)$ of $G$ such that the following holds:
\begin{enumerate}[label={\emph{(\arabic*)}}]
    \item For any $i\in I$ and $j\in J$, either $R_i\subset T_j$ or $R_i\cap T_j = \emptyset$,
    \item $\underline{\mathrm B} (\T ) \geq 1-c-\delta$. 
\end{enumerate}  
\end{lem}

\begin{proof}
Without loss of generality, we may assume $K$ contains the identity element $e$ of $G$ by replacing $K$ with $K \cup \{e\}$ if necessary. By Proposition \ref{prop:disbound}, we can also choose a bounded disjoint quasi-tiling $\Tcl = (J,\Tl,\gamma)$ that is $(\delta, BK)$-boundary-invariant and satisfies   $\underline{\mathrm B} (\T ) \geq 1-\delta$.  For any $j\in J$, define
$$
I_j = \set{i\in I : \alpha_i \in \Tl_j}
\quad \text{and} \quad 
T_j = \p{\Tl_j  \cup \bigcup_{i\in I_j} R_i} \setminus \p{\bigcup_{i\in I \setminus I_j} R_i}.
$$
Set $\T := (J, (T_j), \gamma)$. It is clear that $\T$ is bounded, disjoint, and satisfies Assertion (1).  Moreover together with Lemma \ref{lem:banachunion},  
$$
    \underline{\mathrm B} (\T) \geq 
    \underline{\mathrm B}(\Tcl \cap \R) \geq 
     \underline{\mathrm B} (\mathcal R) + \underline{\mathrm B} (\Tcl ) -1 \geq 1-c-\delta . 
$$
It remains to check the invariance property. By construction we have $\mathrm{Int}_B(\Tl_j) \subset T_j \subset \mathrm{Cl}_B(\Tl_j)$ for any $j\in J$. Then 
\[
T_j K \subset \mathrm{Cl}_B(\Tl_j) K \subset \Tl_j K \cup (\bigcup \{\gamma BK :\gamma B \cap \Tl_j \neq \emptyset , \gamma B \cap \Tl_j^c \neq \emptyset \} ) \subset \mathrm{Cl}_{BK}(\Tl_i) . 
\]
So 
\[ |T_j K| \leq |\mathrm{Cl}_{BK}(\Tl_j) | \leq (1+ \delta) |\mathrm{Int}_{BK}( 
\Tl_j) | \leq  (1+ \delta) |T_j|, 
\]
which implies that $\T$ is $(\delta,K)$-invariant.
\end{proof}

\begin{proof}[Proof of Theorem \ref{thm:wellfiltered}]
	An admissible F\o lner sequence $((F_n)_{n\geq 0}, (B_n)_{n\geq 0}, (\Pc_n)_{n\geq 0})$ can be constructed inductively. 
For $n=0$, we take a bounded partition $\Pc_0$ and let $B_0 \subset G$ be a compact subset such that every atom $Q \in \Pc_0$ is contained in a translation of $B_0$. Let $F_0$ be a sufficiently $B_0$-boundary-invariant F\o lner subset chosen from the given F\o lner sequence.

Let $n \geq 1$, and assume that $F_k, B_k$ and $\Pc_k$ have been well chosen for $0 \leq k < n$, satisfying the required conditions for being a $(1-2^{-n+1})c$-admissible F\o lner sequence, i.e., $(\Pc_n)_{0\leq k<n}$ is a finite atomic filtration and $F_k, B_k$ and $\Pc_k$ satisfy: 
\begin{itemize}
\item for any $0 \leq k < n$ and $A\in \Pc_k$, there exists $\gamma\in G$ such that $A \subset \gamma  B_k$;
\item for any $0 \leq k \leq k' < n$, $F_{k'}$ is $(2^{k-k'},B_k)$-boundary-invariant;
\item  $\underline{\mathrm{B}}(\bigcap_{0 \leq k < n}\bigcup_{T\in\Pc_{k}^{a}}T) \geq 1-(1-2^{-n+1})c$. 
\end{itemize}
To construct $B_n$ and $\Pc_n$, we consider the constant $\delta>0$ and a compact subset $K\subset G$ such that the $(\delta,K$)-invariance is stronger than $(2^{k-n},F_k)$-invariance for any $k < n$. Applying Lemma \ref{lem:tilingforfiltration} to $\Pc_{n-1}$ and $(\delta,K$), we obtain a quasi-tiling $\T_n = (T_i )_{i\in I}$ and a compact subset $B_n \subset G$ such that $\underline{\mathrm{B}}(\T_n ) \geq 1- 2^{-n}c$, and each atom is contained in some translation of $B_n$ with $(\delta,K$)-invariance. Without loss of generality, we may assume that $B_n$ contains a neighborhood of the neutral element $e$ of $G$. Since $G$ is locally compact and second countable, by the same construction as in the proof of Lemma  \ref{lem:tilingforfiltration}, the complement $(\bigcup_{i} T_i )^c$ can be covered by countably many translations of $B_n$ that are nested in $\mathcal P _{n-1}$. So the subset $(\bigcup_{i} T_i )^c$ can be partitioned into atoms $(T_j ')_{j\in J}$, each contained in some translation of $B_n$. In this way, we obtain a bounded partition $\Pc_n = (T_i)\cup (T_j ')$. By Lemma \ref{lem:banachunion}, we have \[ \underline{\mathrm{B}}\left( \bigcap_{k\leq n} \bigcup_{T\in \Pc_k^a}T  \right) \geq \underline{\mathrm{B}}\left( \bigcap_{k\leq n-1} \bigcup_{T\in \Pc_k^a}T  \right) + \underline{\mathrm{B}}\left(  \bigcup_{i \in I}T_i  \right)  -1 =1 - 2^{-n} c .\]Finally, to construct $F_n$, it suffices to choose a sufficiently invariant F\o lner set with respect to $(B_k)_{k\leq n }$. 
\end{proof}

\section{Calderón-Zygmund decomposition}

\label{sec:CZdec}

In this section, we discuss the noncommutative, non-doubling Calder\'{o}n-Zygmund decomposition in the setting of admissible F\o lner sequences. The arguments are adapted from \cite{CCP21} (see also \cite{Cad18,Par09}). As this is one of the only noncommutative extensions of Calder\'on-Zygmund theory in an abstract geometric setting without a doubling condition (or any condition of polynomial growth), we include a self-contained proof for the sake of completeness and clarity. 

To avoid unnecessary analysis of convergence, we consider finite truncation of a given filtered F\o  lner sequence $((F_n)_{n\in I},(B_n)_{n\in I},(\Pc_n)_{n\in I})$, where  $I = \set{0,\dots,N}$ is a finite index set with arbitrarily large $N$. This standard simplification is of no consequence for proving Theorem \ref{thm:maindiffop} since for any sequence $(f_n)_{n\geq 0}$ in $L_{1,\infty}(\N;\ell_2^{rc})$,
$$
\norm{(f_n)_{n\geq 0}}_{L_{1,\infty}(\N;\ell_2^{rc})} =
\sup_{N\geq 0} \norm{(f_n)_{0\leq n\leq N}}_{L_{1,\infty}(\N;\ell_2^{rc})}.
$$ All the constants in the inequalities in this and the next sections are independent of $N$.

Let $(\Eb_n)_{n\in I}$ be the conditional expectations associated with $(\Pc_n)_{n\in I}$. Let $f\in L_{1}(\N)\cap\N$ be a positive and compactly supported function and set $f_n = \mathbb E _n f$.
Let us recall the following construction due to Cuculescu \cite{Cuc71} which is the starting point to define noncommutative Calderón-Zygmund
decompositions.  
\begin{lem}[Cuculescu's projections] \label{lem:cuculescu}
Let $\lambda>0$. The sequence of projections
$(q_{k})_{k\in I}$ determined inductively by $$q_N = \Ind_{[0,\lambda]}(f_N)\quad
\text{and} \quad
q_k = \Ind_{[0,\lambda]}(q_{k+1} f_k q_{k+1}),\quad 0\leq k<N$$ satisfies  
\begin{enumerate}[label={\emph{(\arabic*)}}]
\item for all $0\leq k \leq N$, $q_{k}f_{k}q_{k}\leq\lambda$, 
\item for $k < N$, $q_{k}$ commutes with $q_{k+1}f_{k}q_{k+1}$ and $q_k \leq q_{k+1}$, 
\item $q_{k}\in\N_{k}$, 
\item $\lambda\varphi(1-q_{0})\leq\norm f_{1}.$ 
\end{enumerate}
\end{lem} 
Define, for all $N > n\geq0$, 
\[
p_{n}:=q_{n+1}-q_{n} \quad \text{and} \quad p_N = 1- q_N.
\]
Note that $\sum_{k\geq0}^N p_{k}=1-q_{0}$. We may view $q_k$ or $p_k \in \mathcal N _k $ as an $\mathcal M$-valued function and write 
\[ 
q_k = \sum_{A\in \mathcal P _k } q_A \Ind_A,\quad p_k = \sum_{A\in \mathcal P _k } p_A \Ind_A , \quad q_A, p_A \in \M .
\]
Note that $p_k$ is supported on admissible atoms if $f$ is.

The function $f$ will be decomposed into three parts (plus $p_Nf_N$)
\[
f = g+h+b+p_Nf_N,
\]
a good part, a hybrid part
(which appears when dealing with non-doubling filtrations) and a bad
part. The bad part splits into diagonal terms and off-diagonal terms.
The behaviour of the diagonal terms is similar to the one of the bad part of
the classical Calder\'on-Zygmund decomposition. The off-diagonal terms, however, are
purely noncommutative: they would vanish in a commutative algebra.
Parcet was the first to provide an estimate for those terms in \cite{Par09}
by introducing a pseudo-localisation principle for singular integrals.
In this paper, we follow another strategy presented in 
\cite{CCP21} where off-diagonal terms are grouped
differently and estimated through a more direct computation. The remainder $p_Nf_N$ does not usually appear and it could be removed by making further assumptions on the filtration. In this paper it is more convenient to keep it since it is going to be easy to estimate. 

We collect, in the following lemmas, the construction and essential
properties of each part. \begin{lem}[The good part] \label{lem:good}
The good part, defined by
\begin{equation*}
 g=q_{0}f_{0}q_{0}+\sum_{k = 0}^{N-1} g_{k}  \ \text{ with }\  g_{k}=\Eb_{k+1}(p_{k}f_{k}p_{k}),
 \end{equation*}
satisfies $\norm g_{2}^{2}\leq2\lambda\norm f_{1}$. \end{lem}
\begin{proof} Recall that by Cuculescu's construction (Lemma \ref{lem:cuculescu}),
$q_{k}$ and $q_{k+1}f_{k}q_{k+1}$ commute. Hence, $p_{k}f_{k}p_{k}=q_{k+1}f_{k}q_{k+1}-q_{k}f_{k}q_{k}$
and 
\begin{equation}
g_{k}=q_{k+1}f_{k+1}q_{k+1}-\Eb_{k+1}(q_{k}f_{k}q_{k}).\label{eq:exprforgk}
\end{equation}
We will estimate the $L_2$-norm of $g$ by decomposing $g$ into orthogonal parts, $g = \Eb_N(g) + \sum_{k=0}^{ N-1} \Delta_k(g)$ where $\Delta_k = \Eb_k - \Eb_{k+1}$. We have
\begin{align*}
\Delta_{k}(g) & =\Delta_{k}(q_{0}f_{0}q_{0})+\sum_{i\geq0}\Delta_{k}(g_{i})\\
 & =\Delta_{k}(q_{0}f_{0}q_{0})+\sum_{i=0}^{k-1}\Delta_{k}(q_{i+1}f_{i+1}q_{i+1})-\Delta_{k}(q_{i}f_{i}q_{i})\\
 & =\Delta_{k}(q_{0}f_{0}q_{0})+\Delta_{k}(q_{k}f_{k}q_{k})-\Delta_{k}(q_{0}f_{0}q_{0})=\Delta_{k}(q_{k}f_{k}q_{k})
\end{align*}
and $\Eb_N(g) = q_Nf_Nq_N$ by the same computation. 

Hence,
\[
\norm{\Delta_{k}(g)}^{2}_2=\norm{\Delta_{k}(q_{k}f_{k}q_{k})}_{2}^{2}=\norm{q_{k}f_{k}q_{k}}_{2}^{2}-\norm{\Eb_{k+1}(q_{k}f_{k}q_{k})}_{2}^{2}.
\]
Summing over $k$, we obtain 
\begin{align*}
\norm g_{2}^{2} & = \norm{q_Nf_Nq_N}_2^2 + \sum_{k=1}^{N-1}\norm{q_{k}f_{k}q_{k}}_{2}^{2}-\norm{\Eb_{k+1}(q_{k}f_{k}q_{k})}_{2}^{2}\\
 & \leq \norm{q_{0}f_{0}q_{0}}_{2}^{2}+\sum_{k = 1}^N \left(\norm{q_{k}f_{k}q_{k}}_{2}^{2}-\norm{\Eb_{k}(q_{k-1}f_{k-1}q_{k-1})}_{2}^{2}\right).
\end{align*}
And now let us note that 
\begin{align*}
&\ \quad\norm{q_{k}f_{k}q_{k}}_{2}^{2}-\norm{\Eb_{k}(q_{k-1}f_{k-1}q_{k-1})}_{2}^{2}\\
& =\varphi((q_{k}f_{k}q_{k}-\Eb_{k}(q_{k-1}f_{k-1}q_{k-1}))(q_{k}f_{k}q_{k}+\Eb_{k}(q_{k-1}f_{k-1}q_{k-1})))\\
 & \leq\norm{q_{k}f_{k}q_{k}-\Eb_{k}(q_{k-1}f_{k-1}q_{k-1})}_{1}\norm{q_{k}f_{k}q_{k}+\Eb_{k}(q_{k-1}f_{k-1}q_{k-1})}_{\infty}.
\end{align*}
By Eq.  \eqref{eq:exprforgk}, 
\[
\norm{q_{k}f_{k}q_{k}-\Eb_{k}(q_{k-1}f_{k-1}q_{k-1})}_{1}=\norm{\Eb_{k}(p_{k-1}f_{k-1}p_{k-1})}_{1}=\varphi(p_{k-1}f).
\]
And by the construction of Cuculescu (Lemma \ref{lem:cuculescu}),
\[
\norm{q_{k}f_{k}q_{k}+\Eb_{k}(q_{k-1}f_{k-1}q_{k-1})}_{\infty}\leq2\lambda.
\]
Note also that by the same lemma it holds that $\norm{q_{0}f_{0}q_{0}}_{2}^{2}\leq\varphi(q_{0}f)\lambda$.
Putting the pieces back together we get 
\[
\norm g_{2}^{2}\leq\lambda\varphi(q_{0}f)+\sum_{k\geq1}2\lambda\varphi(p_{k-1}f)\leq2\lambda\norm f_{1},
\]
which establishes the lemma.
\end{proof} \begin{lem}[The hybrid part] \label{lem:hybrid}
The hybrid part, defined by\begin{equation*}
h=\sum_{k\geq 0}^{N-1} h_{k} \ \text{with}\  h_{0}= q_0  (f-f_0) q_0  + p_{0}f_{0}p_{0}-g_{0},\, h_{k}=p_{k}f_{k}p_{k}-g_{k}\text{ for }k\geq 1 ,
\end{equation*}
satisfies 
\begin{enumerate}[label={\emph{(\arabic*)}}]
\item $\Eb_{k+1}(h_{k})=0$, 
\item $h_{k}\in\N_{k}$, 
\item $\sum_{k\geq0}\norm{h_{k}}_{1}\leq2\norm f_{1}$. 
\end{enumerate}
\end{lem} \begin{proof}By the bimodule property of conditional expectations, we see that
\[\mathbb E _{k+1} (q_0 f_0 q_0 ) = \mathbb E _{k+1} (\mathbb E _{0} (q_0 f  q_0) ) = \mathbb E _{k+1} (q_0 f q_0 ). \]By definition $\mathbb E _{k+1} (p_k f_k p_k - g_k )=0$ for $k\geq 0$. Thus the first assertion holds true. 
The second assertion is a straightforward
consequence of the construction. The third is checked as follows.
First note that for any $k\geq0$, 
\[
\norm{p_{k}f_{k}p_{k}- h_{k}}_{1}\leq\norm{p_{k}f_{k}p_{k}}_{1}+\norm{\Eb_{k+1}(p_{k}f_{k}p_{k})}_{1}=2\varphi(p_{k}f).
\]
Then, summing up over $k$ we get 
\[
\sum_{k\geq0}\norm{h_{k}}_{1}\leq2\varphi(q_0 f)+2\sum_{k\geq0}\varphi(p_{k}f) \leq2\varphi( f)=  2\norm f_{1}. \qedhere
\]
\end{proof} \begin{lem}[The bad part] \label{lem:bad} The bad
part, defined by
\begin{equation*}
b=\sum_{k=0}^{N-1}b_{k}^{{\rm d}}+b_{k}^{{\rm off}}\qquad \text{with}\quad  b_{k}^{{\rm d}}=p_{k}(f-f_{k})p_{k}\ \text{ and }\  b_{k}^{{\rm off}}=p_{k}fq_{k}+q_{k}fp_{k} , \end{equation*}
satisfies the following:
\begin{enumerate}[label={\emph{(\arabic*)}}]
\item $\norm{b_k^{\mathrm{d}}}_{1}\leq2\varphi(p_{k}f)$, 
\item let $E$ be a union of $\mathcal{P}_{k}$-atoms and $K\subset E$, 
\begin{equation}
\norm{\int_{K}b_k^{\mathrm{off}}}_{1}\leq\lambda2\varphi(\Ind_{E}p_{k}f).\label{eq:boff}
\end{equation}

\item $\Eb_{k}(b_k^{\mathrm{d}})=\Eb_{k}(b_k^{\mathrm{off}})=0$. 
\end{enumerate}
\end{lem} \begin{proof}The first and third properties follow directly from the definition, hence we focus on the second. By the triangle inequality,
it suffices to consider the case where $E=A$ is an atom of $\mathcal{P}_{k}$.
By definition, $b_k^{\mathrm{off}}=p_{k}fq_{k}+q_{k}fp_{k}$. For convenience,
we only write the argument for $p_{k}fq_{k}$, the other summand admits
the same estimate by adjunction. By the Hölder type inequality stated
in Eq. \eqref{eq:holder}, we have 
\[
\norm{\int_{K}p_{k}fq_{k}}_{L_1 (\mathcal M )}\leq\norm{\p{\int_{K}p_{k}fp_{k}}^{1/2}}_{L_1 (\mathcal M )}\norm{\p{\int_{K}q_{k}fq_{k}}^{1/2}}_{\mathcal M}.
\]
We estimate each factor separately. Note that $\Ind_K p_k = \Ind_K p_A$ by our assumption $E=A$. By  applying  Hölder's inequality once more,
we obtain 
\[
\norm{\p{\int_{K}p_{k}fp_{k}}^{1/2}}_{L_1 (\mathcal M )}\leq\norm{\int_{A}p_{A}fp_{A}}_{L_{1/2} (\mathcal M )}^{1/2}\leq\tau(p_{A})^{1/2}\varphi(\Ind_{A}p_{k}f)^{1/2}.
\]
Given that $q_{A}f q_{A}\leq\lambda$, we can deduce 
\[
\norm{\p{\int_{K}q_{k}fq_{k}}^{1/2}}_{\infty}\leq\norm{\int_{A}q_{k}fq_{k}}_{\infty}^{1/2}\leq\p{\md A\lambda}^{1/2}.
\]
Note that, by Cuculescu's construction (Lemma \ref{lem:cuculescu}), $p_{k}f_{k}p_{k}\geq\lambda p_{k}$, so we have $$\md A\tau(p_{A})\lambda=\varphi(\Ind_{A}p_{k}\lambda)\leq\varphi(\Ind_{A}p_{k}f).$$
Combining the estimates for both terms yields 
\[
\norm{\int_{K}p_{k}fq_{k}}_{1}\leq\tau(p_{A})^{1/2}\varphi(\Ind_{A}p_{k}f)^{1/2}\p{\md A\lambda}^{1/2}\leq\varphi(\Ind_{A}p_{k}f).   \qedhere
\]
\end{proof} 

\begin{lem}[Dilated support] \label{lem:zeta} Assume that $ ((F_n)_{n\geq 0}, (B_n)_{n\geq 0}, ( \mathcal P _n)_{n\geq 0} )$ is an admissible F\o lner sequence and $f$ is also admissible. Set
\[
\zeta=(\bigvee_{k\geq1}\bigvee_{A\in\mathcal{P}_{k}}p_{A}\Ind_{A \p{\cup_{n<k}F_{n}}^{-1}})^{\bot}\in\N.
\]
We have: 
\begin{enumerate}[label={\emph{(\arabic*)}}]
\item $\lambda\varphi(1-\zeta)\les\norm f_{1}$, 
\item for any $x\in G$, $n<k$ and $y\in xF_{n}$,
\[
\zeta(x)b_k^{\mathrm{d}}(y)\zeta(x)=\zeta(x)b_k^{\mathrm{off}}(y)\zeta(x)=0.
\]
\end{enumerate}
\end{lem}
\begin{proof}
	(1) Recall that any admissible atom $A\in \mathcal P _k$ is $(2^{n-k},F_n)$-invariant for $n<k$. We have
 \begin{align*}
 	\varphi(1-\zeta) & \leq \sum_{k\geq 1} \sum_{A\in \mathcal P _k} \tau (p_A ) |A \p{\cup_{n<k}F_{n}}^{-1} | 
 	\leq \sum_{k\geq 1} \sum_{A\in \mathcal P _k} \tau (p_A ) \left(|A|+ \sum_{n<k} |A F_{n}^{-1} \setminus A| \right)\\
 	& \leq 2 \sum_{k\geq 1} \sum_{A\in \mathcal P _k} \tau (p_A ) |A| = 2 \sum_{k\geq 1} \varphi (p_k)
 	=2 \varphi (1-q_0) \\
 	&\leq 2\|f\|_1 / \lambda.
 \end{align*}  
(2) If $y\in A$ for some $A\in \mathcal{P}_k$, then we have $x\in AF_n^{-1}$ and therefore $\zeta (x)  $ is by definition a projection in $\mathcal M$ dominated by $  p_A ^{\bot} $. Note that by definition $$p_A ^{\bot} b_k^{\mathrm{d}}(y)p_A ^{\bot} =p_A ^{\bot} b_k^{\mathrm{off}}(y)p_A ^{\bot} =0.$$
Thus the assertion is verified.
\end{proof}

\section{The difference operator}

\label{sec:diffop}

As in the previous section, we continue to work with an admissible  filtered F\o lner sequence and keep the same notation. In this section, we establish  boundedness properties for the difference operator
\[
\mathbb D:L_1 (\mathcal N ) \to 
 L_{1,\infty} (\mathcal N \overline{\otimes} L_\infty ([0,1])),\quad f\mapsto\sum_{n\in I}\e_{n}\mathbb{D}_{n} f
\]
where $\mathbb{D}_{n}=\mathbb{A}_{n}-\Eb_{n}$ and $(\e_{n})_{n\in I}$ is a sequence
of independent Rademacher variables.  Our analysis is divided into two subsections.  In the first subsection, we derive local estimates, \textit{i.e.}, estimates
for $\mathbb{D}_{n}$ with $n$ fixed. These estimates are primarily based on geometric considerations.
In the second, we prove our main theorem on the weak $(1,1)$ boundedness
of $\mathbb D$ by combining local estimates with the noncommutative Calderón-Zygmund
decomposition introduced earlier.

At its core, our argument is inspired by the approach in \cite{HX21,Xu21}, while we adopt a group-theoretic framework that abstracts away from the Euclidean metric structure. We emphasize the separation between geometric observations and Calderón-Zygmund techniques, and translate the necessary components to make the argument applicable to filtered Følner sequences. Notably, Lemma \ref{lem:wellfilteredL1} is introduced to compensate for the fact that not every atom may be admissible.

\subsection{Local estimates}
In this subsection, the invariance conditions are used to obtain exponential bounds for certain operators.  

\label{sub:local} 
For convenience, given a compact subset $F\subset G$, we write 
$$
[\Ab_F (f)](x) := \dfrac1{\md{F }}f \ast \chi_F(x) = \dfrac1{\md{F }}\int_{F } f(xh)dm(h),\quad f\in L_1(\N)  ,\ x \in G.
$$
Note that $\Ab_{F_n} = \Ab_n $.
\begin{lem}\label{eq:invariance}
If $E$ is $(\e,K)$-invariant then
\begin{equation*}
\norm{\Ind_{E}-\mathbb{A}_{K}(\Ind_{E})}_{1}\leq 2\e\md{E}.
\end{equation*}
\end{lem}

\begin{proof} For any $y\in K$, we have
$$\md{Ey\Delta E} = \md{Ey \setminus E} +  \md{Ey^{-1} \setminus E} \leq 2\e\md{E}.$$
Then desired inequality follows from the computation: 
\begin{align*}
\norm{\Ind_{E}-\mathbb{A}_{K}(\Ind_E)}_{1} & =\int_{G}\md{\Ind_{E}(x)-\md K^{-1}\int_{K}\Ind_{E}(xy)dm(y)}dm(x)\\
 & =\md K^{-1}\int_{G}\md{\int_{K}(\Ind_{E}(x)-\Ind_{E}(xy))dm(y)}dm(x)\\
 & \leq\md K^{-1}\int_{G}\int_{K}\md{\Ind_{E}(x)-\Ind_{E}(xy)}dm(y)dm(x)\\
 & =\md K^{-1}\int_{K}\md{E\Delta Ey}dm(y)\\
 & \leq\md{2K} ^{-1}\int_{K}\e\md Edm(y)=2\e\md E.    \qedhere
\end{align*}
\end{proof} 

\begin{prop} \label{prop:localestimates} 
Let $p\in[1,\infty)$ and $n,k\geq0$. Let  $f\in L_{p}(\N)$ be an admissible function.
\begin{enumerate}[label={\emph{(\arabic*)}}]
\item if $n<k$ and $f\in L_{p}(\N_{k})$, then $\norm{\mathbb{D}_{n}(f)}_{p}\les2^{(n-k)/p}\norm f_{p}$, 
\item if $n\geq k$ and $\Eb_{k}(f)=0$, then $\norm{\mathbb{D}_{n}(f)}_{p}\les2^{k-n}\norm f_{p}$. 
\end{enumerate}
\end{prop} 

\begin{proof} (1) Note that $\norm{\mathbb{D}_{n}:\N_{k}\to\N_{k}}\leq 2$.
So by complex interpolation and considering the subalgebra of functions supported on admissible atoms $\mathcal P _k^a$ in $\N_k$, it suffices to treat the case $f\in L_{1}(\N_{k})$.
In this case we may assume that $f=f_{A}\Ind_{A}$ for an atom $A\in\mathcal{P}_{k}$
and that $f_{A}$ is positive. 
Note that $A $ belongs to the $\sigma$-subalgebra generated by $\mathcal P_n$ with $n<k$, which in particular
implies that $\mathbf{1}_{A}\in\mathcal{N}_{n}$ and hence 
\[
\mathbb{D}_{n}(f)=f_{A}(\mathbb{A}_{n}(\mathbf{1}_{A})-\mathbb{E}_{n}(\mathbf{1}_{A}))=f_{A}(\mathbb{A}_{n}(\mathbf{1}_{A})-\Ind_{A}).
\]
Using Lemma \ref{eq:invariance} and the $(2^{n-k},F_{n})$-invariance,
we get 
\[
\norm{\mathbb{D}_{n}(f)}_{1}=\norm{f_{A}}_{1}\norm{\mathbb{A}_{n}(\mathbf{1}_{A})-\Ind_{A}}_{1}\les2^{n-k}\norm f_{1}.
\]

(2) Assume now that $n\geq k$. Note that $\Eb_{k}(f)=0$, so $\int_{A}f=0$
for all $A\in\mathcal{P}_{k}$ and $\mathbb{E}_{n}(f)=\mathbb{E}_{n}\mathbb{E}_{k}(f)=0$.
Therefore for $x\in G$,

\[
\mathbb{D}_{n}(f)(x)=\mathbb{A}_{n}(f)(x)=\dfrac{1}{\md{xF_{n}}}\int_{xF_{n}}f=\dfrac{1}{\md{xF_{n}}}\int_{\partial_{\mathcal{P}_{k}}(xF_{n})\cap xF_{n}}f\eqqcolon \mathbb{D}_{n,k}(f)(x),
\]
where $\partial_{\mathcal{P}_{k} }(xF_{n}):=\bigcup\set{A\in\mathcal{P}_{k} :A\cap xF_{n}\neq\emptyset,A\nsubseteq xF_{n}}$.
Recall that by definition of an admissible sequence, any
$A\in\mathcal{P}_{k} $ is contained in a left-translate of $B_{k}$ and
therefore $\partial_{\mathcal{P}_{k} }(xF_{n})\subset\partial_{B_{k}}(xF_{n})$.

Let us show that $\norm{\mathbb{D}_{n,k}:L_{p}(\N)\to L_{p}(\N)}\leq2^{k-n}$.
Since $\mathbb{D}_{n,k}$ is positive, let us consider $\rho\in L_{p}(\N)^+$. We have for $x\in G$,
\begin{align*}
\mathbb{D}_{n,k}(\rho)(x) & \leq\dfrac{1}{\md{xF_{n}}}\int_{\partial_{B_{k}}(xF_{n})}\rho=\dfrac{\md{\partial_{B_{k}}(xF_{n})}}{\md{xF_{n}}}\mathbb{A}_{\partial_{B_{k}}(F_{n})}(\rho)(x).
\end{align*}
Note that $\partial_{B_{k}}(xF_{n})=x\partial_{B_{k}}(F_{n})$ and
$F_{n}$ is by assumption $(2^{k-n},B_{k})$-boundary-invariant. We therefore
obtain
\[
\mathbb{D}_{n,k}(\rho)\leq\dfrac{\md{\partial_{B_{k}}(F_{n})}}{\md{F_{n}}}\mathbb{A}_{\partial_{B_{k}}(F_{n})}(\rho)\leq2^{k-n}\mathbb{A}_{\partial_{B_{k}}(F_{n})}(\rho).
\]
Since $\mathbb{A}_{\partial_{B_{k}}(F_{n})}$ is a contraction on $L_{p}(\N)$
for all $p\in[1,\infty]$, this concludes the proof. \end{proof}
Let us collect two more local estimates which are consequences of
specific properties of the Calderón-Zygmund decomposition presented
in the previous section. Fix an admissible function $f\in L_{1}(\N)^+ \cap \N$ with compact support. 

\begin{lem} \label{lem:localcancellation}
For two integers $n<k$, 
\[
\zeta \mathbb{D}_{n}(b_k^{\mathrm{d}})\zeta=\zeta \mathbb{D}_{n}(b_k^{\mathrm{off}})\zeta=0.
\]
\end{lem} \begin{proof} Note that by Lemma \ref{lem:zeta}, for
$x\in G$, 
\[
\zeta(x)\mathbb{A}_{n}(b_k^{\mathrm{d}})(x)\zeta(x)=\md{F_{n}}^{-1}\int_{F_{n}}\zeta(x)b_k^{\mathrm{d}}(xy)\zeta(x)dm(y)=0.
\]
Moreover, by the definition
of $\zeta$ in Lemma \ref{lem:zeta}, we see that $\zeta^{\bot}\geq p_{k}$.
But $p_{k}\in\mathcal{N}_{k}\subset\mathcal{N}_{n}$ since $n<k$.
So together with the definition of $b_k^{\mathrm{d}}$, we have $\Eb_{n}(b_k^{\mathrm{d}})=\Eb_{n}(p_{k}b_k^{\mathrm{d}})=p_{k}\Eb_{n}(b_k^{\mathrm{d}})$.
Hence, 
\[
\zeta\Eb_{n}(b_k^{\mathrm{d}})\zeta=\zeta p_{k}\Eb_{n}(b_k^{\mathrm{d}})\zeta=0.
\]
The same argument applies to $b_k^{\mathrm{off}}$. 
\end{proof} 

\begin{lem}\label{lem:localboff} 
For $n\geq k$, $\norm{\mathbb{D}_{n}(b_k^{\mathrm{off}})}_{1}\les2^{k-n}\varphi(p_{k}f)$.
\end{lem} 

\begin{proof} Since $\Eb_{n}(b_k^{\mathrm{off}})=0$ (see Lemma
\ref{lem:bad}), using the notation and observations made in the
proof of Proposition \ref{prop:localestimates}, we have
\[
\mathbb{D}_{n}(b_k^{\mathrm{off}})(x)=\md{xF_{n}}^{-1}\int_{xF_{n}\cap\partial_{\mathcal{P}_{k} }(xF_{n})}b_k^{\mathrm{off}},
\]
Then, by Eq. \eqref{eq:boff} in Lemma \ref{lem:bad} applied to $K = xF_{n}\cap\partial_{\mathcal{P}_{k} }(xF_{n})$ and $E = \partial_{\mathcal{P}_{k} }(xF_{n})$
\begin{align*}
\norm{\mathbb{D}_{n}(b_k^{\mathrm{off}})(x)}_{1} & \leq\md{xF_{n}}^{-1}\tau\p{\int_{\partial_{\mathcal{P}_{k} }(xF_{n})}p_{k}fp_{k}}.
\end{align*}
Recall that as in the proof of Proposition \ref{prop:localestimates},
we have $\partial_{\mathcal{P}_{k} }(xF_{n})\subset\partial_{B_{k}}(xF_{n})$,
$\partial_{B_{k}}(xF_{n})=x\partial_{B_{k}}(F_{n})$ and that $F_{n}$
is $(2^{k-n},B_{k})$-boundary-invariant. Thus the above inequality yields
\[
\norm{\mathbb{D}_{n}(b_k^{\mathrm{off}})}_{1}
\leq \dfrac{|\partial_{B_{k}}(F_{n})|}{\md{F_{n}}}\norm{\mathbb{A}_{\partial_{B_{k}}(F_{n})}(p_{k}f)}_{1}
\leq 2^{k-n}\varphi(p_{k}f).    \qedhere
\]
\end{proof}

\subsection{Weak $(1,1)$ boundedness}

Now we are ready to prove the desired boundedness property of the
difference operator $\mathbb{D}$. We start by establishing the $L_{2}$-boundedness
of $\mathbb{D}$. \begin{prop} \label{prop:L2boundforD} Let $f\in L_{2}(\N)$ be an admissible function. Then 
\[
\norm{\mathbb{D}(f)}_{2}\les\norm f_{2}.
\]
\end{prop} \begin{proof}Recall that it suffices to obtain uniform bounds on the finite sum $\sum_{n = 0}^N \varepsilon_n \mathbb D _n$. Decompose
$$f = f_N + \sum_{k=-1}^{N-1} df_{k}  \quad \text{with }  df_{k}=f_{k}-f_{k+1},\ f_{-1}\coloneqq f.
$$
The $L_2$ estimate for $\Db(f_N)$ can be directly derived from Proposition \ref{prop:localestimates} so we focus on $f' := \sum_{k=-1}^{N-1} df_{k}$. Note that $\Eb_{k+1}(df_{k})=0$
and $df_{k}\in\N_{k}$ so that by the two local estimates stated in
Proposition \ref{prop:localestimates}, for any $k \in \set{-1,\dots,N-1}$ and $n\in\set{0,\dots,N}$, we have
\[
\norm{\mathbb{D}_{n}(df_{k})}_{2}^{2}\lesssim2^{-\md{k-n}}\norm{df_{k}}_{2}^{2}
\]
Set $df_k = 0$ for $k\notin \set{-1,\N-1}$. We have
\begin{align*}
\norm{\mathbb{D}(f')}_{2} 
& =\left(\sum_{n= 0}^N\|\mathbb{D}_{n}f'\|_{2}^{2}\right)^{1/2}
=\left(\sum_{n=0}^N \norm{\sum_{k=-1}^{N-1}\mathbb{D}_n(df_{k})}_{2}^{2}\right)^{1/2}\\
 & =\left(\sum_{n=0}^N\|\sum_{s \in \Zb} \mathbb{D}_{n}(df_{n+s})\|_{2}^{2}\right)^{1/2}
 \les\left(\sum_{n=0}^N\left(\sum_{s\in\Zb}2^{-\md{s}/2}\|df_{n+s}\|_{2}\right)^{2}\right)^{1/2}.
\end{align*}
Hence by Young's inequality for convolutions, 
\[
\norm{\mathbb{D}(f')}_{2} \les \p{\sum_{s\in\Zb} 2^{-\md{s}/2}}\p{\sum_{k=-1}^{N-1}\norm{df_k}_2^2}^{1/2} \les \norm{f'}_2.   \qedhere
\]
\end{proof}

We are nearly in a position to prove Theorem \ref{thm:maindiffop}. The following lemma helps to reduce the discussion of Theorem  \ref{thm:maindiffop} to the case of admissible functions. 
\begin{lem} \label{lem:wellfilteredL1} Let $E\subset G$
be a compact subset with a nonzero Banach density. Let $f\in L_{1}(\N)$ be of compact support
(as a function from $G$ to $L_{1}(\M)$). Then there exists a sequence
$(f_{n})_{n\in\Nb}$ in $L_{1}(\N)$, a sequence $(x _{n})_{n\geq0}\subset G$ and a constant $\alpha>0$ such that 
\begin{enumerate}[label={\emph{(\arabic*)}}]
\item for any $n\geq0$, $f_{n}$ is supported in $x _n E$, 
\item for any $n\geq0$, $\norm{f_{n}}_{1}\leq \alpha^{-n}\norm f_{1}$, 
\item $\sum_{n=0}^{\infty}f_{n }=f$. 
\end{enumerate}
\end{lem} 

\begin{proof}
Let $K$ be the support of $f$. Assume without loss of generality that $f$ is positive. For any $\e >0$, we can find an $(\e,K^{-1})$-invariant
compact subset $D$ of $G$. We may also assume, by the nonzero Banach density condition, that $\md{E \cap D} \geq c \md{D}$ for some $c \in (0,1)$. In particular, for any $y \in K$,
$$
\md{Ey^{-1} \cap D  K^{-1}} \geq \md{(E\cap D)y^{-1}} \geq c \md{D}. 
$$
We aim to find an element $x_0$ of $G$ such that $\norm{f\Ind_{x_0^{-1}E}}_1 \geq c'\norm{f}_1 $ for some $c'\in(0,1)$. To do so, we compute the average value of $\norm{f\Ind_{x^{-1}E}}_1$ for $x \in D  K^{-1}$
\begin{align*}
\frac1{\md{D  K^{-1}}}\int_{x\in D  K^{-1}}\norm{f\Ind_{x^{-1}E}}_1 dr(x) & =\dfrac1{\md{D  K^{-1}}}\tau \p{\int_{y\in K}f(y)\int_{x\in D  K^{-1}}\Ind_{Ey^{-1}}(x)dr(x)dr(y)}\\
 &= \dfrac1{\md{D  K^{-1}}}\int_{y\in K} \tau(f(y)) \md{D  K^{-1} \cap Ey^{-1}} dr(y)\\
 &\geq \dfrac{c\md D}{ \md{D  K^{-1}}}\norm f_{1} \geq \dfrac{ c}{ 1+\e} \norm{f}_1
\end{align*}
which implies there exists $x_{0}\in D  K^{-1}$ such that 
\[
\norm{f\Ind_{x_0^{-1}E}}_1 \geq \dfrac{ c}{ 1+\e} \norm f_{1}.
\]
Set $f_{0}=f\Ind_{x_{0}^{-1}E}$ and $f'=f-f_{0}$. Note that $\norm{f'}_{1}\leq \frac{1+\e -c}{1+\e }\norm f_{1}$.
We can then repeat the argument, substituting $f'$ for $f$ to construct $f_{1}$.
The result follows by induction.
\end{proof}

\begin{proof}[Proof of Theorem \ref{thm:maindiffop}] Let us show that we may assume $f$ to be admissible. First, by density, $f$ can be taken to be of compact support. Then by Lemma \ref{lem:wellfilteredL1}, decompose $f=\sum_n f_n$ for $E = \mathrm{Adm}(\mathcal S )$ with $\norm{f_{n}}_{1}\leq \alpha^{-n}\norm f_{1}$ for all $n$. In particular each function $f_n$ is admissible. Assume that the theorem holds for admissible functions. Then  for each $n$ with $\lambda_{n}=\alpha^{-n/2}\lambda$,
we obtain projections $e_{n}$ satisfying
\begin{equation*}
 \norm{\left(\sum_k \varepsilon_k \mathbb D _k f_n \right)e_{n}}_{\infty}\leq \alpha^{-n}\lambda 
\quad{}\text{and} \quad{} \lambda\varphi(1-e_{n})\leq \alpha^{-n/2}\norm f_{1}.\end{equation*}
Set $e=\bigwedge_{n\geq0}e_{n}$. Then we see that 
\begin{equation*}
 \norm{\left(\sum_k \varepsilon_k\mathbb D _k f    \right)e }_{\infty}\lesssim  \lambda 
\quad{}\text{and} \quad{} \lambda\varphi(1-e )\lesssim \norm f_{1} .\end{equation*}
In other words we have $ \| \mathbb D f \|_{L_{1,\infty} (\mathcal N ;\ell_2^{rc})}\lesssim \|f\|_1$. 

Thus in the following we consider the case where $f\in L_{1}(\N)^{+}$ is a bounded and compactly
supported admissible function. Replace without loss of generality $\Db$ by its truncated version $\Db = \sum_{n=0}^N \e_n\Db_n$. 

Let $\lambda>0$. We aim to show that 
\[
\lambda\varphi(\set{\md{\mathbb{D}(f)}>\lambda})\les\norm f_{1}.
\]
By Eq. \eqref{eq:lambda}, it suffices to deal with each part of the Calderón-Zygmund
decomposition separately, namely to show that 
\[
\lambda\varphi(\set{\md{\mathbb{D}(u)}>\lambda})\les\norm f_{1}
\]
for $u=p_Nf_N,g,h$ and $b$. Note that $u$ is also admissible by construction.

\noindent \textit{Estimate for $p_Nf_N$.} By Proposition \ref{prop:localestimates} (1), since $p_nf_N \in \N_N$, 
$$
\norm{\Db(p_Nf_N)} \les \norm{p_Nf_N}_1 \leq \norm{f}_1
$$
and we can conclude using Tchebychev's inequality. 

\noindent \textit{The good part.} By Lemma \ref{lem:good}, $\norm g_{2}^{2}\les\lambda\norm f_{1}$.
Using Tchebychev's inequality as well as the $L_{2}$-boundedness
of $\mathbb{D}$ (Proposition \ref{prop:L2boundforD}) we get 
\[
\varphi(\set{\md{\mathbb{D}(g)}>\lambda})\les\dfrac{\norm{\mathbb{D}(g)}_{2}^{2}}{\lambda^{2}}\les\dfrac{\norm g_{2}^{2}}{\lambda^{2}}\les\dfrac{\norm f_{1}}{\lambda}.
\]

\noindent \textit{The hybrid part.} By Lemma \ref{lem:hybrid}, we
have $\Eb_{k+1}(h_{k})=0$ and $h_{k}\in\N_{k}$. So by the local estimates
in Proposition \ref{prop:localestimates}, for any $k,n\geq0$ 
\[
\norm{\mathbb{D}_{n}(h_{k})}_{1}\les2^{-\md{k-n}}\norm{h_{k}}_{1}.
\]
Summing up those inequalities over $n$ and $k$, we get 
\[
\norm{\mathbb{D}(h)}_{1}\les\sum_{k\geq0}\norm{h_{k}}_{1}.
\]
It follows, using the norm estimate in Lemma \ref{lem:hybrid}, that
\[
\norm{\mathbb{D}(h)}_{1}\les\norm f_{1}.
\]
We conclude once again using Tchebychev's inequality.

\noindent \textit{The bad part.} For the bad part, it is not possible
to get norm estimates right away. However, Lemma \ref{lem:localcancellation}
tells us that $\mathbb{D}$ (and in general Calder\'on-Zygmund-operators) behaves better on the projection $\zeta$ \textit{i.e.} away from the support of
$b$. 

Concretely, by Eq. \eqref{eq:lambda} we may write 
\begin{align*}
\varphi(\set{\md{\mathbb{D}(b)}>\lambda}) & \leq\varphi(\set{\md{\zeta^{\bot}\mathbb{D}(b)}>0})+\varphi(\set{\md{\zeta \mathbb{D}(b)\zeta^{\bot}}>0})+\varphi(\set{\md{\zeta \mathbb{D}(b)\zeta}>\lambda})\\
 & \les\lambda^{-1}\norm f_{1}+\lambda^{-1}\norm{\zeta \mathbb{D}(b)\zeta}_{1}
\end{align*}
where we used Lemma \ref{lem:zeta} to estimate the first two terms
and Tchebychev's inequality for the last one. Hence, it remains to
prove that 
\[
\norm{\zeta \mathbb{D}(b)\zeta}_{1}\les\norm f_{1}.
\]
For the diagonal part, by Lemma \ref{lem:bad}, Lemma \ref{lem:localcancellation}
and Proposition \ref{prop:localestimates} (2), we have $\zeta \mathbb{D}_{n}(b_k^{\mathrm{d}})\zeta =0$ for any $n<k$ and
\[
\norm{\zeta \mathbb{D}_{n}(b_k^{\mathrm{d}})\zeta}_{1}\les2^{-\md{n-k}}\norm{b_k^{\mathrm{d}}}_{1}\lesssim 2^{-\md{n-k}}\lambda\varphi(p_{k})\]for any $n\geq k$. 
Similarly, by Lemma \ref{lem:localboff} and Lemma \ref{lem:localcancellation},
\[
\norm{\zeta\mathbb{D}_{n}(b_k^{\mathrm{off}}) \zeta}_{1}\les2^{-\md{n-k}}\lambda\varphi(p_{k}).
\]
Now note that by Lemma \ref{lem:cuculescu}, 
\[
\lambda\sum_{k\geq0}\varphi(p_{k} )=\lambda\varphi(1-q_{0})\leq\norm f_{1},
\]
which establishes the desired inequality. \end{proof}
 
\subsection*{Acknowledgment}
The authors would like to thank Romain Tessera for several helpful discussions. They would also like to thank Patrick Poissel and Wei Liu for their careful reading of the preprint and for pointing out several mistakes and misprints. The authors were partially supported by public grants as part of the Fondation
Mathématique Jacques Hadamard. L. Cadilhac was also supported by the French (Agence Nationale de la Recherche) grant ANR-19-CE40-0002. S. Wang was also partially
supported by the Fundamental Research Funds for the Central Universities No.
FRFCUAUGA5710012222 and NSF of China No. 12031004.

\appendix
\section{Completely admissible F\o lner sequences}
In this appendix, we show that the construction of admissible F\o lner sequences presented in Section \ref{sec:geometry} can be strengthened if the group $G$ is unimodular. The key observation is the following extension of Proposition \ref{prop:disbound}, where the tiling $\T$ can be indeed chosen to be a partition in the unimodular setting. Part of the proof is adapted from \cite{DHZ19}.

\begin{prop}
    Let $G$ be a unimodular locally compact second countable amenable group. Let $\delta >0$ and $K \subset G$ be a compact subset. There exists a tiling $\T$ of $G$ that is a partition, bounded, and $(\delta,K)$-boundary-invariant.
\end{prop}
\begin{proof}
Let $\e$ be a small positive constant,	and fix  a  $(\varepsilon,K)$-invariant compact subset $T\subset G$ containing the neutral element $\{e\}$. Set $C\coloneqq 8\varepsilon |T|^{-1} $, and let $L=L(T,K)$ be a sufficiently large compact subset of $G$ depending only on $T$ and $K$.
 Let $\T = (I,(T_i)_{i\in I},(\gamma_i)_{i\in I})$ be a disjoint, $B$-bounded, $(\e /2,L)$-boundary-invariant tiling of $G$ with $\bd{\T} \geq 1-\e/3$ given by Proposition \ref{prop:disbound}. Let $F$ be a compact set such that
	\[ \inf_{g\in G}\frac{|\bigcup_i T_i  \cap g  F |}{|g  F|}>1-\frac{\varepsilon}{2}, \]
	and such that $F  $ is $(\e /2,B )$-boundary-invariant and $(\e /2,T)$-invariant. Let 
	\[ E\coloneqq \bigcup_{j\in J} \gamma_j T \subset \T^c \]
	be a maximal disjoint union of translates of $T$ in $\T^c $. Then by maximality, any translate $g  T$  intersects either $\bigcup_i T_i $ or $E$; in other words, $G$ is covered by the union of $\T \cdot T^{-1}$ and $E  T^{-1}$. We may rewrite these two sets as disjoint unions, setting for $i\in I$ and $j\in J$,
	$$
 A_i = T_i  T^{-1} \setminus \bigcup_{i'<i} T_{i'}  T^{-1}
 \quad \text{and} \quad 
 B_j = \gamma_j  T  T^{-1} \setminus \p{\bigcup_{j'<j}\gamma_{j'}  T  T^{-1} \cup \bigcup_i T_i  T^{-1}}.
 $$
Note that for each $i$ and $j$,
	\[ T_i \subset A_i \subset  T_i T^{-1} ,\quad
	\gamma_j T \subset B_j \subset 
	\gamma_j T T^{-1}, \]
which in particular implies that $A_i$ is $(\e/2 ,K)$-invariant by the choice of $T_i$  and $L$.	We claim that there is a map $f:J\to I$ such that 
	\begin{equation}\label{claima1}
		\forall j\in J, \quad \gamma_j  F \cap T_{f(j)} \neq \emptyset\quad \text{ and }\quad
		\forall i\in I ,\quad | f^{-1}(i) |\leq C |T_i |.
	\end{equation} 
To prove this, let us introduce a family $(c_i)_{i\in I}$ such that $c_i  T_i \subset B$. and consider the relation $\mathscr R$ between $J$ and ${ \tilde{{I}}\coloneqq \{ (i,n)\mid n\leq C  |T_i | \}}$ defined as follows: for any $j\in J$ and $(i,n)\in \tilde I$,
	\[ j\mathscr R (i,n) \text{ iff }  c_i \in \gamma_j   F .\] 
	We can assume that $\min_i C|T_i | \geq 1$ (for example by going back to the construction of $\T$ that starts with Proposition \ref{prop:fromOW} and guarantees that we can choose sufficiently large atoms).  

We have the following observations. 
	\begin{itemize}
		\item[a)] For any fixed $j\in J$, we consider
		\begin{align*}
			\left|\{(i,n)\in\tilde I\mid j\mathscr R (i,n)\}\right|& 
			= \sum_{i: c_i \in \gamma_j  F }\lfloor C  |T_i |\rfloor \geq \frac{1}{2}  C \sum_{i: c_i \in \gamma_j    F } |T_i |\\
			&\geq  \frac{1}{2}  C \left(|(\cup_i T_i)\cap \gamma_j   F |- |\partial_B (\gamma_j   F) |\right)\\
			&
			\geq  \frac{1-\varepsilon }{4}  C |  F | ,
		\end{align*}
where the last inequality is due to the choice of $F$.
		\item[b)] For any fixed $(i,n)\in \tilde I$, 
		\begin{align*}
			|\{j\mid j\mathscr R (i,n)\} |
			&=|\{j\mid c_i^{-1}\gamma_j     \in F\}| = |T|^{-1} \sum_{j:c_i^{-1}\gamma_j     \in F}   | \gamma_j T|\\
			&	\leq |T|^{-1} |E\cap   c_i F T | \\
			& \leq   |T|^{-1}( |E\cap   c_i F| +|  c_i  FT \setminus   c_i  F| )\\
			&\leq \varepsilon  |T|^{-1} |F|  .
		\end{align*}
	\end{itemize}
With the above observations a) and b), we may apply the Hall marriage lemma to obtain an injective map from $\tilde I$ to $J$, which asserts \eqref{claima1}.
	
	We now consider the atoms $T_i ' \coloneqq A_i \cup (\bigcup_{j\in f^{-1}(i)}B_j)$. Together with the invariant conditions on $A_i$ and $T$, we see that
	\begin{align*}
		|T_i ' K |&\leq |A_i K | + \sum_{j\in f^{-1}(i)} |TK|
		\leq |A_i K | +   | f^{-1}(i) | |TK| \\
		&\leq |A_i K | +  C |TK| |T_i | \\
        & \leq (1+\frac{\varepsilon}{2}) |A_i | + C (1+\varepsilon ) |T| |T_i | \leq (1+17 \varepsilon ) |T_i ' |.
	\end{align*}
	Thus $T_i '$ is $(\delta,K)$-invariant as soon as $\varepsilon$ was chosen small enough. The proof is complete.
\end{proof}
Then a verbatim but simpler arguments as in the proofs of Lemma \ref{lem:tilingforfiltration} and Theorem \ref{thm:wellfiltered} yields the following result. Following a similar terminology, we say that an  admissible F\o lner sequence $ ((F_n)_{n\geq 0},(B_n)_{n\geq 0},(\mathcal{P}_n)_{n\geq 0})$ is \emph{completely admissible} if every atom of $\mathcal{P}_n$ for any $n\geq 0$ is admissible.

\begin{thm}\label{thm:cfiltered}
Let $G$ be a unimodular locally compact second countable amenable group. Any F\o lner sequence in $G$ admits
a completely admissible F\o lner subsequence. 
\end{thm}

\bibliographystyle{abbrv}
\bibliography{Bibli}

\end{document}